\documentclass[11pt]{article}

\usepackage{amsmath,amssymb,amsthm,amscd}
\usepackage{latexsym}
\usepackage[all]{xy}

\allowdisplaybreaks[1]

 \theoremstyle{plain}
\newtheorem{thm}{Theorem}[section]
\newtheorem{lemma}[thm]{Lemma}
\newtheorem{proposition}[thm]{Proposition}
\newtheorem{cor}[thm]{Corollary}

\theoremstyle{definition}

\xyoption{dvips}
\numberwithin{equation}{section}
\newcommand{\FF}{\mathbb{F}}  
\newcommand{\ZZ}{\mathbb{Z}}  
\newcommand{\CC}{\mathbb{C}}

\newcommand{\GL}{\mathrm{GL}}

\newcommand{\cP}{\mathcal{P}}
\newcommand{\cS}{\mathcal{S}}
\newcommand{\Res}{\mathrm{Res}}
\newcommand{\Inf}{\mathrm{Inf}}
\newcommand{\Sinf}{\mathrm{Sinf}}
\newcommand{\SInd}{\mathrm{SInd}}
\newcommand{\NCSym}{\mathrm{NCSym}}
\newcommand{\scs}{\scriptstyle}

\newcommand{\fkn}{\mathfrak{n}}
\newcommand{\ch}{\mathrm{ch}}
\newcommand{\cK}{\mathcal{K}}
\newcommand{\cX}{\mathcal{X}}
\newcommand{\cC}{\mathcal{C}}

\newcommand{\larc}[1]{\overset{#1}{\frown}}

\newcommand{\One}{{1\hspace{-.14cm} 1}}
\newcommand{\spanning}{\text{-span}}
\newcommand{\dd}{\displaystyle}

\begin{document}

\title{Branching rules in the ring of superclass functions\\ of unipotent upper-triangular matrices}
\author{Nathaniel Thiem\\ Department of Mathematics\\ University of Colorado at Boulder\\
 \textsf{thiemn@colorado.edu}}

\date{}

\maketitle

\begin{abstract}
 It is becoming increasingly clear that the supercharacter theory of the finite group of unipotent upper-triangular matrices has a rich combinatorial structure built on set-partitions that is analogous to the partition combinatorics of the classical representation theory of the symmetric group.  This paper begins by exploring a connection to the ring of symmetric functions in non-commuting variables that mirrors the symmetric group's relationship with the ring of symmetric functions.    It then also investigates some of the representation theoretic structure constants arising from the  restriction, tensor products and superinduction of supercharacters in this context.  
 \end{abstract}

\section{Introduction}


The representation theory of the symmetric group $S_n$ --with its connections to partition and tableaux combinatorics-- has become a fundamental model in combinatorial representation theory.  It has become clear in recent years that the representation theory  finite group of unipotent upper-triangular matrix groups $U_n(q)$ can lead to a similarly rich combinatorial theory.  While understanding the usual representation theory of $U_n(q)$ is a wild problem, Andr\'e \cite{An95,An99,An01,An02} and Yan \cite{Ya01,Ya06} constructed  a natural approximation to the representation theory that leads to a more computable theory.  This approximation (known as a super-representation theory) now relies on set-partition combinatorics in the same way that the representation theory of the symmetric group relies on partition combinatorics.  

A fundamental tool in symmetric group combinatorics is the ring of symmetric functions, which encodes the character theory of all symmetric groups simultaneously in a way that polynomial multiplication in the ring of symmetric functions  becomes  symmetric group induction from Young subgroups.  This kind of a relationship has been extended to wreath products and type $A$ finite groups of Lie type (for descriptions see for example \cite{Ma95,TV07}).   One of the purposes of this paper is to suggest an analogous relationship between the supercharacter theory of $U_n(q)$ and the ring on symmetric functions in non-commuting variables $\NCSym$.  In particular, Corollary \ref{CharacteristicMap} shows that  there are a family of algebra isomorphisms from the ring of supercharacters to $\NCSym$, where we replace induction from subgroups with its natural analogue superinduction from subgroups.  Unfortunately, there does not yet seem to be a canonical choice (ideally, such a choice would take the Hopf structure of $\NCSym$ into account).

The other purpose of this paper is to use the combinatorics of set partitions to supply recursive algorithms for computing restrictions to a family of subgroups called parabolic subgroups.  It turns out that if $k\leq n$, then there are many ways in which $U_k(q)$ sits inside $U_n(q)$ as a subgroup.  In fact, for every subset $S\subseteq \{1,2,\ldots, n\}$ with $k$ elements, there is a distinct subgroup $U_S$ of $U_n(q)$ isomorphic to $U_k(q)$. The restriction from $U_n(q)$ to $U_S$ depends on $S$, and Theorem \ref{RestrictionRule} sorts out the combinatorial differences for all possible subsets $S$.   This result can then be easily extended to give restriction rules for all parabolic subgroups.  These computations require knowledge of tensor product results  that  were previously done by Andr\'e \cite{An95} for large prime and by Yan \cite{Ya01} for arbitrary primes.  For completeness, this paper supplies an alternate proof that relates tensor products to restriction and a generalization of the inflation functor (see Lemma \ref{TensorLemma}).

By Frobenius reciprocity we then also obtain the coefficients of superinduction from these subgroups.    Corollary \ref{SuperInductionTwoParts} concludes by observing that superinduced supercharacters from parabolic subgroups are essentially twisted super-permutation characters (again using the generalization of the inflation map).  These results give the structure constants for the ring of superclass functions of the finite unipotent upper-triangular groups.   However, the underlying coefficient ring is $\ZZ[q^{-1}]$, unlike in the case of the symmetric group where the ring is $\ZZ$.

The paper is organized as follows.  Section \ref{Preliminaries} introduces some set-partition combinatorics; describes the parabolic subgroups that will replace Young subgroups in our theory; reviews the supercharacter theory of pattern groups (as defined in \cite{DI08}); and recalls the ring of symmetric functions in non-commuting variables $\NCSym$.  We proceed in Section \ref{CharacteristicMap} by describing the family of isormorphisms between $\NCSym$ and the ring of supercharacters.  Section \ref{SectionBranchingRules} uses the fact that supercharacters of $U_n(q)$ decompose into tensor products of simpler characters to supply algorithms for computing restrictions and superinductions of supercharacters.  These results generalize restriction results in \cite{TV07}, and make use of a new generalization of the inflation functor to supercharacters of pattern groups.  

This paper builds on \cite{MT07}  and \cite{TV07p} by giving restriction and superinduction formulas for larger families of groups.  These formulas are computable, and are being implemented in Python as part of an honors thesis at the University of Colorado.  Other recent work in this area worth mentioning includes extensions by Andr\'e and his collaborators to supercharacter theories of other types \cite{AN06} and over other rings \cite{AnNi06}, explorations of all supercharacter theories for a given group by Hendrickson in his thesis \cite{He08}, and an intriguing unexplored connection to $L$-packets in the work of Drinfeld and Boyarchenko \cite{BD06}.

\vspace{.25cm}
 
 \noindent\textbf{Acknowledgements.}  Much of this work was completed at MSRI during the simultaneous programs ``Combinatorial Representation Theory" and ``Representation Theory of Finite Groups" in Spring 2008.  I especially appreciated related conversations with Andr\'e, Diaconis, Isaacs, and Yan during this period.  

%
%

\section{Preliminaries}\label{Preliminaries}

This section reviews the combinatorics needed for the main results, gives a brief introduction to the supercharacter theory of pattern groups, and recalls the ring of symmetric functions in non-commuting variables.

%
%

\subsection{$\FF_q$-labeled set-partitions}\label{SetPartitionCombinatorics}

For $S\subseteq \{1,2,\ldots, n\}$, let 
$$\cS_S  = \{\text{set-partitions of $S$}\},$$
and
$$\cS=\bigcup_{n\geq 0} \cS_n,\qquad\text{where}\qquad \cS_n  = \cS_{\{1,2,\ldots, n\}}.$$

An \emph{arc} $i\larc{}j$ of $K\in \cS_S$ is a pair $(i,j)\in S\times S$ such that 
\begin{enumerate}
\item[(1)] $i<j$,
\item[(2)] $i$ and $j$ are in the same part of $K$,
\item[(3)] if $k$ is in the same part as $i$ and $i<k\leq j$, then $k=j$.
\end{enumerate}
Thus, if we order each part in increasing order, then the arcs are pairs of adjacent elements in each part.  For example,
$$\{1,5,7\}\cup\{2,3\}\cup \{4\} \cup \{6,8,9\}\in \cS_9$$
has arcs $1\larc{}5$, $5\larc{}7$, $2\larc{}3$, $6\larc{} 8$, and $8\larc{}9$.  We can also represent the set partition $K$ as a diagram consisting of $|K|$ vertices with an edge connecting vertex $i$ to vertex $j$ if $i\larc{}j$ is an arc of $K$; for example,
$$\{1,5,7\}\cup\{2,3\}\cup \{4\} \cup \{6,8,9\}
\qquad\longleftrightarrow\qquad 
\xymatrix@R=.2cm@C=.5cm{*={\bullet} \ar @{-} @/^{.5cm}/ [rrrr] & *={\bullet}  \ar @{-} @/^{.2cm}/ [r]  & *={\bullet} & *={\bullet} & *={\bullet} \ar @{-} @/^{.3cm}/ [rr] & *={\bullet}  \ar @{-} @/^{.3cm}/ [rr]  & *={\bullet} & *={\bullet}  \ar @{-} @/^{.2cm}/ [r] & *={\bullet}\\
*={\scs 1} & *={\scs 2} & *={\scs  3} & *={\scs 4} & *={\scs 5} & *={\scs 6} & *={\scs 7} & *={\scs 8} & *={\scs 9} }\ .$$
The \emph{arc set} $A(K)$ of $K\in \cS_S$ is
$$A(K)=\{\text{arcs of $K$}\}.$$

A \emph{crossing} of $K\in \cS_S$ is a pair of arcs $(i\larc{}k,j\larc{}l)\in A(K)\times A(K)$ such that $i<j<k<l$.   The \emph{crossing set} $C(K)$ of $K$ is
$$C(K)=\{\text{crossings of $K$}\}.$$
For example, if $K=\{1,5,7\}\cup\{2,3\}\cup \{4\} \cup \{6,8,9\}$, then $K$ has one crossing $(5\larc{} 7, 6\larc{} 8)$, as is easily observed in the above diagrammatic representation of $K$.

An \emph{$\FF_q$-labeled set-partition of $S$} is a pair $(\lambda,\tau_\lambda)$, where $\lambda$ is a set-partition of $S$ and $\tau_\lambda:A(\lambda)\rightarrow \FF_q^\times$ is a labeling of the arcs by nonzero elements of $\FF_q$.  By convention, if $\tau_\lambda(i\larc{}j)=a$, then we write the arc as $i\larc{a} j$.  Thus, we can typically suppress the labeling function in the notation.  Let
$$\cS_S(q)=\{\text{$\FF_q$-labeled set-partitions of $S$}\},$$
and
$$\cS(q)=\bigcup_{n\geq 0} \cS_n(q), \qquad \text{where}\quad \cS_n(q)=\cS_{\{1,2,\ldots,n\}}(q).$$ 
Note that if $s_n(q)=|\cS_n(q)|$, then the generating function
$$\sum_{n\geq 0} s_n(q) \frac{x^n}{n!}=e^{\frac{e^{(q-1)x}-1}{q-1}}$$
is a $q$-analogue of the usual exponential generating function of the Bell numbers (where $q=2$ gives the usual generating function).

Suppose $S\subseteq T\subseteq \{1,2,\ldots,n\}$.  Then there is a function
$$\begin{array}{rccc} \langle\cdot\rangle_T: & \left\{\begin{array}{c} \text{$\FF_q$-labeled}\\ \text{set-partitions of $S$}\end{array}\right\}   &\longrightarrow  &\left\{\begin{array}{c} \text{$\FF_q$-labeled}\\ \text{set-partitions of $T$}\end{array}\right\}\\
& \lambda & \mapsto & \langle \lambda\rangle_T\end{array}$$
where $\langle \lambda\rangle_T$ is the unique $\FF_q$-labeled set-partition of $T$ with arc set $A(\lambda)$ and labeling function $\tau_\lambda$.   We will use the convention that $\langle \lambda\rangle_n=\langle \lambda\rangle_{\{1,2,\ldots, n\}}$.

%
%

\subsection{Pattern groups}

For $n\in \ZZ_{\geq 1}$, let $U_n(q)$ be the group of $n\times n$ unipotent upper-triangular matrices with entries in $\FF_q$.  Given a poset $\cP$ of $\{1,2,\ldots, n\}$, the \emph{pattern group} $U_\cP(q)$ is 
$$U_\cP(q)=\{u\in U_n(q) \mid u_{ij}\neq 0 \text{ implies $i\preceq j$ in } \cP\}.$$

\noindent\textbf{Remark.}  If $T_n(q)$ is the group of $n\times n$ diagonal matrices with entries in $\FF_q^\times$, then the set of pattern subgroups of $U_n(q)$ can be characterized as the set of subgroups fixed by the conjugation action of $T_n(q)$ on $U_n(q)$. 

\vspace{.25cm}

Consider the injective map 
$$\begin{array}{ccc} \cS_n & \longrightarrow & \left\{\begin{array}{c} \text{Posets of}\\ \{1,2,\ldots, n\}\end{array}\right\}\\
K & \mapsto & \cP_K\end{array}$$
where $i\prec j$ in $\cP_K$ if and only if $i<j$ and both $i$ and $j$ are in the same part of $K$.

A pattern subgroup $U_{\cP}(q)$  is a \emph{parabolic subgroup} of $U_n(q)$ if there exists $K\in \cS_n$ such that $\cP=\cP_K$.  Note that if $K=K_1\cup K_2\cup\cdots \cup K_\ell$ is the decomposition of $K$ into parts, then 
$$U_{\cP_K}(q)\cong U_{|K_1|}(q)\times U_{|K_2|}(q)\times \cdots \times U_{|K_\ell|}(q).$$
Thus, the parabolic subgroups of $U_\cP(q)$ are reminiscent of the Young subgroups of the symmetric groups $S_n$ or parabolic subgroups of the general linear groups $\GL_n(q)$. In fact, we will follow this analogy into the supercharacter theory of $U_n(q)$.  To simplify notation, we will typically write
$$U_K(q)=U_{\cP_K}(q), \qquad \text{for $K\in \cS_n$.}$$

\vspace{.25cm}

\noindent\textbf{Remarks.}
\begin{enumerate}
\item[(a)]   These subgroups are not generally block diagonal.  For example, 
$$U_{\cP_{\{1,3,5\}\cup \{2,4\}}} =\left\{\left(\begin{array}{ccccc} 
1 & 0 & \ast & 0 & \ast\\ 
0 & 1 & 0 & \ast & 0\\
0 & 0 & 1 & 0 & \ast\\
0 & 0 & 0 & 1 & 0 \\
0 & 0 & 0 & 0 & 1\end{array}\right) \bigg| \ast\in \FF_q\right\}\cong U_3(q)\times U_2(q).$$
\item[(b)]  Parabolic subgroups do not include all possible copies of pattern subgroups isomorphic to a direct product of $U_k(q)$'s.  For example,
$$U_{\xymatrix@R=.2cm@C=.2cm{*={} & *{\scs 4}\\ *{\scs 3}  \ar @{-} [ur] & *={} & *{\scs 2}\ar @{-} [ul]\ar @{-} [dl]\\ *={} & *{\scs 1}}}=\left\{\left(\begin{array}{cccc} 
1 & 0 & \ast &  \ast\\ 
0 & 1 & 0 & \ast\\
0 & 0 & 1 &  \ast\\
0 & 0 & 0 & 1 \end{array}\right) \bigg| \ast\in \FF_q\right\}\cong U_3(q)\times U_2(q)$$
is not a parabolic subgroup of $U_4(q)$.
\end{enumerate}

%
%

\subsection{A supercharacter theory for pattern groups}

Given a group $G$, a supercharacter theory is an approximation to the usual character theory.  To be more precise, a supercharacter theory consists  of a set of superclasses $\cK$ and a set of supercharacters $\cX$, such that
\begin{enumerate}
\item[(a)] the set $\cK$ is a partition of $G$ such that each part is a union of conjugacy classes,
\item[(b)] the set $\cX$ is a set of characters such that each irreducible character appears as the constituent of exactly one supercharacter,
\item[(c)] the supercharacters are constant on superclasses,
\item[(d)] $|\cK|=|\cX|$,
\item[(e)] the identity element $1$ of $G$ is in its own superclass, and the trivial character $\One$ of $G$ is a supercharacter. 
\end{enumerate} 
This general notion of a supercharacter theory was introduced by Diaconis and Isaacs \cite{DI08} to generalize work of Andr\'e and Yan on the character theory of $U_n(q)$.  

\vspace{.25cm}

\noindent\textbf{Remark.}  The definition includes a reasonable amount of redundancy, as explored in \cite{DI08,He08}.

\vspace{.25cm}

Diaconis and Isaacs extended the construction of Andr\'e of a supercharacter theory for $U_n(q)$ to a larger family of groups called algebra groups.  We will review the construction for pattern groups (a subset of the set of algebra groups).  Let $\cP$ be a poset of $\{1,2,\ldots, n\}$ and let
$$\fkn_\cP(q)=U_\cP(q)-1,$$
which is an $\FF_q$-algebra.  

Fix a nontrivial homomorphism $\vartheta:\FF_q^+\rightarrow \CC^\times$.   The pattern group $U_\cP(q)$ acts on the left and right on both $\fkn_\cP(q)$ and the dual space $\fkn_\cP(q)^\ast$, and the two-sided orbits lead to the sets $\cK$ and $\cX$ by the following rules.
The superclasses are given
\begin{equation*}
\begin{array}{ccc} U_\cP(q)\backslash \fkn_\cP(q)/U_\cP(q) & \longleftrightarrow & \cK\\ U_\cP(q) X U_\cP(q) & \mapsto & 1+U_\cP(q) X U_\cP(q),\end{array}
\end{equation*}
and the supercharacters are given by
$$\begin{array}{ccc} U_\cP(q)\backslash \fkn_\cP(q)^*/U_\cP(q) & \longleftrightarrow & \cX\\ U_\cP(q) \lambda U_\cP(q) & \mapsto & \dd \chi^\lambda=\frac{|\lambda U_\cP(q)|}{|U_\cP(q)\lambda U_\cP(q)|}\sum_{\mu\in U_\cP(q)\lambda U_\cP(q)}\vartheta\circ \mu. \end{array}
$$
The corresponding $U_\cP$-modules are given by
$$V^\lambda=\CC\spanning\{v_\mu\mid \mu\in U_\cP \lambda\},$$
with action
$$gv_\mu=\vartheta\big((g\mu)(1-g)\big)v_{g\mu}, \qquad \text{for $g\in U_\cP$ and $\mu\in U_\cP \lambda$.}$$

\vspace{.25cm}

\noindent\textbf{Examples.}  The group $U_n(q)$ was the original motivation for studying supercharacter theories.  The following results are due to Andr\'e, Yan,  and Arias-Castro--Diaconis--Stanley.  The number of superclasses is
$$|\cK|=|\cX|=|\cS_n(q)|,$$
where, for example,
$$\begin{array}{ccc} \cS_n(q) & \longrightarrow & \cK\\ \mu & \mapsto & u_\mu,\end{array}\qquad\text{and}\quad (u_\mu)_{ij}=\left\{\begin{array}{ll} 1, & \text{if $i=j$,}\\ \tau_\mu(i\larc{}j), & \text{if $i\larc{}j\in A(\mu)$,}\\ 0, & \text{otherwise.}\end{array}\right.$$
The corresponding supercharacter formula for $\lambda,\mu\in \cS_n(q)$ is
\begin{equation} \label{CharacterFormula}
\chi^{\lambda}(u_\mu)=\left\{\begin{array}{ll}\dd\prod_{i\larc{}l\in A(\lambda)} \frac{q^{l-i-1}\vartheta\big(\tau_\lambda(i\larc{}l)\tau_\mu(i\larc{}l)\big)}{q^{|\{j\larc{}k\in A(\mu)\mid i<j<k<l\}|}},  & \begin{array}{@{}l}\text{if $i<j<k$, $i\larc{}k\in A(\lambda)$}\\ \text{implies $i\larc{}j,j\larc{}k\notin A(\mu)$,}\end{array}\\ 0, & \text{otherwise,}\end{array}\right.
\end{equation}
where $\tau_\mu(i\larc{}j)=0$ if $i\larc{}l\notin A(\mu)$ (see \cite{DT07} for the corresponding formula for arbitrary pattern groups).  Note that the degree of each character is
\begin{equation} \label{DegreeFormula}
\chi^\lambda(1)=\prod_{i\larc{}l\in A(\lambda)} q^{l-i-1}.
\end{equation}
It follows directly from the formula that the supercharacters factor nicely
$$\chi^\lambda=\prod_{i\larc{a}l\in A(\lambda)} \chi^{\langle i\larc{a} l\rangle_n}.$$
It also follows from (\ref{CharacterFormula}) and (\ref{DegreeFormula}) that $\chi^\lambda$ is linear if and only if
$$i\larc{}k\in A(\lambda) \qquad\text{implies}\qquad k=i+1.$$

The set $C(\lambda)$ measures how close the supercharacter $\chi^\lambda$ is to being irreducible.  In fact,
\begin{equation}\label{OrthogonalityRelation}
\langle \chi^\lambda, \chi^\mu\rangle=q^{|C(\lambda)|} \delta_{\lambda\mu},
\end{equation}
where $\langle\cdot,\cdot\rangle$ is the usual inner product on characters.

For parabolics subgroups $U_{K}(q)$ of $U_n(q)$, 
$$|\cK|=|\cX|= |\cS_{|K_1|}(q)||\cS_{|K_2|}(q)|\cdots |\cS_{|K_\ell|}(q)|,$$
where $K=K_1\cup K_2\cup\cdots\cup K_\ell\in \cS_n$.

\vspace{.25cm}

\noindent\textbf{Remark.}  If instead of considering $U_n(q)$-orbits on $\fkn_n(q)$ and $\fkn_n(q)^*$, we consider orbits of the full Borel subgroup $B_n(q)=T_n(q)U_n(q)$ on these spaces, then the corresponding supercharacter theory no longer depends on the finite field $q$.  In this case, the combinatorics reduces to considering set-partitions rather than $\FF_q$-labeled set-partitions.

\vspace{.25cm}

Supercharacters satisfy a variety of nice properties, as described in \cite{DI08}.  The above construction satisfies 
\begin{enumerate}
\item[(a)] The product of two supercharacters is a $\ZZ_{\geq 0}$-linear combination of supercharacters.
\item[(b)] The restriction of a supercharacter from one pattern group to a pattern subgroup is a $\ZZ_{\geq 0}$-linear combination of supercharacters.
\end{enumerate}
However, it is not true that the induction functor sends a supercharacter to a  $\ZZ_{\geq 0}$-linear combination of supercharacters.  In fact, an induced supercharacter is generally no longer even a superclass function.  

Diaconis and Isaacs therefore define a map superinduction on supercharacters that is adjoint to restriction with respect to the usual inner product on class functions; it turns out that this function averages over superclasses in the same way induction averages over conjugacy classes.  In particular, if $H\subseteq G$ are pattern groups (or more generally algebra groups), then \emph{superinduction} is the function
$$\begin{array}{rccc} \SInd: & \left\{\begin{array}{c} \text{Superclass functions}\\ \text{of $H$}\end{array}\right\} & \longrightarrow &  \left\{\begin{array}{c} \text{Superclass functions}\\ \text{of $G$}\end{array}\right\}\\ & \chi & \mapsto & \SInd_H^G(\chi),\end{array}$$
where
$$\SInd_H^G(\chi)(g)=\frac{1}{|G||H|}\sum_{x,y\in G\atop x(g-1)y+1\in H} \chi(1+x(g-1)y), \qquad \text{for $g\in G$.}$$
Unfortunately, while $\SInd$ sends superclass functions to superclass functions, it sends supercharacters to $\ZZ_{\geq 0}[1/q]$-linear combinations of supercharacters (where $q$ comes from the underlying finite field).  In fact, the image is not even generally a character.   See also \cite{MT07} for a further exploration of the relationship between superinduction and induction.

%
%

\subsection{The ring of symmetric functions in non-commutative variables}

Fix a set $X=\{X_1,X_2,\ldots\}$ of countably many non-commuting variables.  For $K=K_1\cup K_2\cup\cdots\cup K_\ell\in \cS_n$, define the \emph{monomial symmetric function} 
$$m_K(X)=\sum_{k=(k_1,k_2,\ldots, k_\ell)\in \ZZ_{\geq 1}^\ell\atop k_i\neq k_j, 1\leq i<j\leq \ell} X_{\pi_1(k)} X_{\pi_2(k)}\cdots X_{\pi_\ell(k)},\qquad\text{where $\pi_j(k)=k_i$ if $j\in K_i$.}$$
The \emph{space of symmetric functions in non-commuting variables of homogeneous degree $n$} is
$$\NCSym_n(X)=\CC\spanning\{m_K(X) \mid K\in \cS_n\},$$ 
and the \emph{ring of symmetric functions in non-commuting variables} is
$$\NCSym=\bigoplus_{n\geq 0} \NCSym_n(X),$$
where a possible multiplication is given by usual polynomial products.   However, note that if $K=\{a_1<a_2<\cdots< a_m\}\cup \{b_{1}<b_{2}<\cdots< b_{n}\}\in \cS_{m+n}$ with $w=(a_{1},a_{2},\cdots,a_{k_m},b_1,b_2,\ldots, b_n)$ the corresponding permutation of $m+1$ elements, then we could ``shuffle" two words according to $K$,
$$(X_{i_1}X_{i_2}\cdots X_{i_m}) \ast_K (X_{i_{m+1}}\cdots X_{i_{m+n}})= X_{i_{w^{-1}(1)}}X_{i_{w^{-1}(2)}}\cdots X_{i_{w^{-1}(m+n)}}.$$
These operations give a variety of alternate shuffle products for $\NCSym$.

The ring $\NCSym$ naturally generalizes the usual ring of symmetric functions \cite{Ma95}, but is different from other generalizations such as the ring of noncommutative symmetric functions studied in, for example, \cite{GKLLRT95}.  The ring $\NCSym$ was introduced by Wolf \cite{Wo36}, and further explored by Rosas and Sagan \cite{RS04}.  There has been recent interest in the Hopf structure of $\NCSym$ and its Hopf dual -- for example, \cite{BHRZ06, BRRZ08}.  In particular,   \cite{BHRZ06} show that it has a representation theoretic connection with partition lattice algebras.    This paper suggests that the supercharacter theory of $U_n(q)$ also has a representation theoretic connection to $\NCSym$ in a way that is more analogous to how the ring of symmetric functions dictates the representation theory of $S_n$.   However, the precise nature of this connection remains open.  In particular, it is not clear whether the Hopf structure of $\NCSym$ translates naturally into a representation theoretic Hopf structure for the supercharacters of $U_n(q)$.

%
%

\section{The ring of unipotent superclass functions}\label{CharacteristicMap}

This section explores the relationship between $\NCSym$ and the space of supercharacters
$$\cC(q)=\bigoplus_{n\geq 0} \cC_n(q), \qquad\text{where}\quad \cC_n(q)=\CC\spanning\{\chi^\lambda\mid\lambda\in \cS_n(q)\}.$$

%
%

\subsection{Parabolic subgroups and set-partition combinatorics}

There are different copies of $U_m(q)\times U_n(q)$ as subgroups of $U_{m+n}(q)$ which are not related via an inner automorphism of $U_{m+n}(q)$.  In fact, for every $K=K_1\cup K_2\in \cS_{m+n}$ with $|K_1|=m$ and $|K_2|=n$, $U_{m+n}(q)$ has a parabolic subgroup $U_m(q)\times_K U_n(q)=U_{K}(q)\cong U_{m}(q)\times U_n(q)$.

Thus, the space $\cC$ has a variety of different products.  For $\lambda\in \cS_m(q)$, $\mu\in \cS_n(q)$, and $K=K_1\cup K_2\in \cS_{m+n}$ with $|K_1|=m$ and $|K_2|=n$, define
$$\chi^\lambda\ast_K \chi^\mu=\SInd_{U_m(q)\times_K U_n(q)}^{U_{m+n}(q)}(\chi^\lambda\times \chi^\mu).$$
There is a related map 
$$\begin{array}{rccc} \cup_K: & \cS_m(q)\times\cS_n(q) & \longrightarrow & \cS_{m+n}(q)\\ 
& (\lambda,\mu) & \mapsto & \lambda \cup_K \mu,\end{array}$$
where $\lambda\cup_K\mu=\lambda'\cup \mu'$ with $\lambda'\in \cS_{K_1}(q)$ and $\mu'\in \cS_{K_2}(q)$ the same $\FF_q$-labeled set-partitions as $\lambda$ and $\mu$ respectively, but with $\{1,2,\ldots, m\}$ relabeled as $K_1$ and $\{1,2,\ldots, n\}$ relabeled as $K_2$.  For example, 
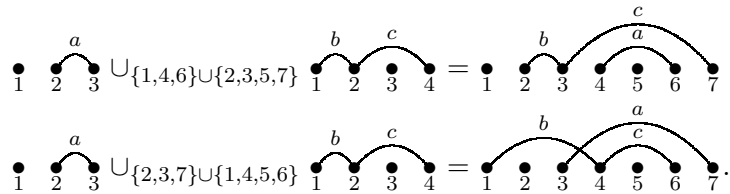
\begin{align*}
\xymatrix@R=.2cm@C=.5cm{*={\bullet} & *={\bullet}  \ar @{-} @/^{.2cm}/ [r]^{a}  & *={\bullet} \\
*={\scs 1} & *={\scs 2} & *={\scs  3}} 
\ &\cup_{\{1,4,6\}\cup\{2,3,5,7\}}\ 
\xymatrix@R=.2cm@C=.5cm{*={\bullet} \ar @{-} @/^{.2cm}/ [r]^{b} & *={\bullet}  \ar @{-} @/^{.3cm}/ [rr]^{c}  & *={\bullet} & *={\bullet} \\
*={\scs 1} & *={\scs 2} & *={\scs  3} & *={\scs 4} }
\ =\ 
\xymatrix@R=.2cm@C=.5cm{*={\bullet} & *={\bullet}  \ar @{-} @/^{.2cm}/ [r]^{b}  & *={\bullet} \ar @{-} @/^{.6cm}/ [rrrr]^{c}  & *={\bullet} \ar @{-} @/^{.3cm}/ [rr]^{a}  & *={\bullet} & *={\bullet}   & *={\bullet} \\
*={\scs 1} & *={\scs 2} & *={\scs  3} & *={\scs 4} & *={\scs 5} & *={\scs 6} & *={\scs 7}  }
\\
\xymatrix@R=.2cm@C=.5cm{*={\bullet} & *={\bullet}  \ar @{-} @/^{.2cm}/ [r]^{a}  & *={\bullet} \\
*={\scs 1} & *={\scs 2} & *={\scs  3}} 
\ &\cup_{\{2,3,7\}\cup\{1,4,5,6\}}\ 
\xymatrix@R=.2cm@C=.5cm{*={\bullet} \ar @{-} @/^{.2cm}/ [r]^{b} & *={\bullet}  \ar @{-} @/^{.3cm}/ [rr]^{c}  & *={\bullet} & *={\bullet} \\
*={\scs 1} & *={\scs 2} & *={\scs  3} & *={\scs 4} }
\ =\ 
\xymatrix@R=.2cm@C=.5cm{*={\bullet} \ar @{-} @/^{.4cm}/ [rrr]^{b} & *={\bullet}   & *={\bullet} \ar @{-} @/^{.6cm}/ [rrrr]^{a}  & *={\bullet} \ar @{-} @/^{.3cm}/ [rr]^{c}   & *={\bullet} & *={\bullet}   & *={\bullet} \\
*={\scs 1} & *={\scs 2} & *={\scs  3} & *={\scs 4} & *={\scs 5} & *={\scs 6} & *={\scs 7}  }\ .
\end{align*}
It will follow from Corollary \ref{SuperInductionTwoParts} that $\chi^{\lambda\cup_K \mu}$ is always a nonzero constituent of $\chi^\lambda\ast_K\chi^\mu$.  

\vspace{.25cm}

\noindent\textbf{Remark.}  The graph automorphism $\sigma$ of the Dynkin diagram of type $A$ gives a natural map
$$\begin{array}{rcccl}\sigma: & \cC & \longrightarrow & \cC\\ 
& \chi^\lambda & \mapsto & \chi^{\sigma(\lambda)}, & \text{for $\lambda\in \cS(q)$,} \end{array}
$$
where $\sigma(\lambda)\in \cS(q)$ is the $\FF_q$-labeled set partition obtained by reflecting the diagram $\lambda$ across a vertical axis.  This map is an anti-automorphism of $\cC$ (it also sends left modules to right modules).   Thus, in many of the following results we only prove half of the symmetric cases.

%
%

\subsection{A characteristic map for supercharacters}

For $\mu\in \cS_n(q)$, let $\kappa_\mu:U_n\rightarrow U_n$ be the superclass characteristic function given by 
$$\kappa_\mu(u)=\left\{\begin{array}{ll} 1, & \text{if $u$ is in the same superclass as $u_\mu$,}\\ 0, & \text{otherwise,}\end{array}\right.$$
and
$$z_\mu=\frac{|U_n|}{|U_n(u_\mu-1) U_n|}$$

\begin{proposition}
For $\mu\in \cS_m(q)$ and $\nu\in \cS_n(q)$, 
$$\SInd_{U_m\times_K U_n}^{U_{m+n}} \big( (z_\mu \kappa_\mu)\otimes(z_\nu \kappa_\nu)\big)=z_{\mu\cup_K\nu} \kappa_{\mu\cup_K\nu}.$$
\end{proposition}
\begin{proof}
By definition,
\begin{align*}
\SInd_{U_m\times_K U_n}^{U_{m+n}} \big( (z_\mu \kappa_\mu)\otimes(z_\nu \kappa_\nu)\big)(g) & = \frac{z_\mu z_\nu}{|U_{m+n}||U_m||U_n|}\sum_{x,y\in U_{m+n}\atop x(g-1)y+1\in U_m\times_K U_n}\hspace{-.5cm} \vartheta(\kappa_\mu\otimes \kappa_\nu)(x(g-1)y+1)\\
&=0,
\end{align*}
unless the superclass containing $g$ also contains $u_{\mu\cup_K\nu}=u_\mu\times_K u_\nu$.  That is, there exists $c\in \CC$ such that
$$\SInd_{U_m\times_K U_n}^{U_{m+n}} \big( (z_\mu \kappa_\mu)\otimes(z_\nu \kappa_\nu)\big)=c \kappa_{\mu\cup_K\nu}.$$
Specifically,
\begin{align*}
c & = \frac{z_\mu z_\nu}{|U_{m+n}||U_m||U_n|}\sum_{x,y\in U_{m+n}\atop x(g-1)y+1\in U_m\times_K U_n} \vartheta(\kappa_\mu\otimes \kappa_\nu)(x(u_{\mu\cup_K \nu}-1)y+1)\\
& = \frac{z_\mu z_\nu}{|U_{m+n}||U_m||U_n|}|U_{m+n}|z_{\mu\cup_K\nu}\sum_{g-1\in U_{m+n}(u_{\mu\cup_K\nu}-1)U_{m+n}\atop g\in U_{m}\times_K U_n} \vartheta(\kappa_\mu\otimes \kappa_\nu)(g)\\
& = \frac{z_\mu z_\nu}{|U_{m+n}||U_m||U_n|}|U_{m+n}|z_{\mu\cup_K\nu}|U_m(u_\mu-1)U_m||U_n(u_\nu-1)U_n|\\ 
&=z_{\mu\cup_K\nu},
\end{align*}
as desired.
\end{proof}

Let $\NCSym$ be the ring of symmetric functions in non-commuting variables.  Let 
$$\{p_\lambda \mid \lambda\in \cS\}$$
be any basis that satisfies
$$p_\lambda \ast_K p_\mu=p_{\lambda\cup_K \mu}$$
for all $K=K_1\cup K_2\in \cS$ with $|K_1|=|\lambda|$ and $|K_2|=|\mu|$.
Note that $\NCSym$ in fact has several of such bases, such as $\{p_\lambda\}$ in \cite{RS04} and $\{x_\lambda\}$ in \cite{BHRZ06}.

\begin{cor} \label{CharacteristicMap}
The function
$$\begin{array}{rccc} \ch: & \cC(2) & \longrightarrow & \NCSym\\ & \kappa_\mu & \mapsto & \frac{1}{z_\mu}p_\mu\end{array}$$
is an isometric algebra isomorphism.
\end{cor}

\noindent\textbf{Questions.}  This result raises the following questions.
\begin{enumerate}
\item[(1)] Does the Hopf algebra structure of $\NCSym$ transfer in a representation theoretic way to $\cC$?
\item[(2)] What is the correct choice of basis $p_\mu$?  In particular, the $\{p_\lambda\}$ of \cite{RS04} do not seem to give a nice Hopf structure to $\cC$.
\item[(3)] Is there a corresponding $\NCSym$-space for $q>2$?
\end{enumerate}
Questions (1) and (2)  presumably need simultaneous answers, and question (3) suggests  there might be an  analogue of the ring symmetric functions corresponding to wreath products.

%
%

\section{Representation theoretic structure constants}\label{SectionBranchingRules}

This section explores the computation of structure constants in $\cC$.   We begin with a family of natural embedding  maps of $\cC_m(q)\subseteq \cC_n(q)$ for $m\leq n$ using a generalization of the inflation functor, and then give algorithms for computing restrictions from $\cC_{m+n}(q)$ to $\cC_m(q)\otimes \cC_n(q)$.  To finish the computations we require a method for decomposing tensor products $\cC_n(q)\otimes \cC_n(q)\rightarrow \cC_n(q)$.  We conclude with a discussion of the corresponding superinduction coefficients.  In this section we will assume a fixed $q$, and suppress the $q$ from the notation; that is $U_n=U_n(q)$, etc.

%
%

\subsection{Superinflation of characters}

Let $T\subseteq G$ be pattern groups with corresponding algebras $\mathfrak{t}$ and $\mathfrak{g}$, respectively.   There exists a surjective projection
$$\begin{array}{rccc} \pi: & \mathfrak{g}=\mathfrak{t}\oplus \mathfrak{t}^\perp & \longrightarrow & \mathfrak{t}\\
&X+Y & \mapsto & X,\end{array}$$
with a corresponding inflation map
$$\begin{array}{rccc} \Inf_{\mathfrak{t}}^{\mathfrak{g}}: & \mathfrak{t}^* & \longrightarrow & \mathfrak{g}^*\\  &  \mu & \mapsto & \mu\circ\pi.\end{array}$$
The \emph{superinflation} map on supermodules is given by
$$\begin{array}{rccc} \Sinf_T^G: & \left\{\begin{array}{c} \text{Supermodules}\\ \text{of $T$}\end{array}\right\} & \longrightarrow & \left\{\begin{array}{c} \text{Supermodules}\\ \text{of $G$}\end{array}\right\}\\
& V^\mu & \mapsto & V^{\Inf_{\mathfrak{t}}^{\mathfrak{g}}(\mu)}.\end{array}$$

Note that superinflation takes supermodules to supermodules, just as the usual inflation map on characters takes irreducible characters to irreducible characters.  Recall, the usual inflation map is constructed from a surjection $\pi:G\rightarrow T$ is given by
$$\begin{array}{rccc} \Inf_T^G: & \{\text{$T$-modules}\} & \longrightarrow & \{\text{$G$-modules}\} \\  & V & \mapsto & \Inf_T^G(V),\end{array}$$
where $gv=\pi(g) v$ for $g\in G$, $v\in  \Inf_T^G(V)$.   The following proposition says that superinflation is inflation whenever possible.

\begin{proposition}
Suppose $G$ is a pattern group with pattern subgroups $T$ and $H$ such that $G=T\ltimes H$.  Then for any supermodule $V^\lambda$ of $T$,
$$\Sinf_T^G(V^\lambda)\cong\Inf_T^G(V^\lambda).$$ 
\end{proposition}

\begin{proof}  Let $\mathfrak{g}=G-1$, $\mathfrak{h}=H-1$, and $\mathfrak{t}=T-1$.  Consider the map
$$\begin{array}{rccc}\varphi: & V^{\Inf_{\mathfrak{t}}^{\mathfrak{g}}(\lambda)} & \longrightarrow & \Inf_T^G(V^\lambda)\\ 
& v_\mu & \mapsto & v_{\Res_{\mathfrak{t}}^{\mathfrak{g}}(\mu)}\end{array}$$ 
Since $T\lambda\subseteq G\lambda$, this map is surjective. 

By \cite[Lemma 3.2]{MT07} normality in pattern groups implies ``super-normality" in the sense that for $h\in H$ and $g\in G$
$$g(h-1), (h-1)g\in \mathfrak{h}.$$
Thus, for $t\in T$ and $h\in H$, 
$$\pi(th-1)=\pi(t(h-1)+(t-1))=\pi(t(h-1))+\pi(t-1)=t-1,$$
and similarly $\pi(ht-1)=t-1$.
For $s,t\in T$, and $h,k\in H$,
\begin{align*}
(ht)\Inf_{\mathfrak{t}}^{\mathfrak{g}}(\lambda)(ks-1)&=\Inf_{\mathfrak{t}}^{\mathfrak{g}}(\lambda)(t^{-1}h^{-1}ks-t^{-1}h^{-1})\\
&=\Inf_{\mathfrak{t}}^{\mathfrak{g}}(\lambda)(t^{-1}s s^{-1}h^{-1}ks-1)+\Inf_{\mathfrak{t}}^{\mathfrak{g}}(\lambda)(1-t^{-1}h^{-1})\\
\intertext{and since $s^{-1}h^{-1}ks\in H$,}
(ht)\Inf_{\mathfrak{t}}^{\mathfrak{g}}(\lambda)(ks-1)&=\lambda(t^{-1}s-1)+\lambda(1-t^{-1})\\
&=(t\lambda)(s-1)\\
&=t\Res_{\mathfrak{t}}^{\mathfrak{g}}\big(\Inf_{\mathfrak{t}}^{\mathfrak{g}}(\lambda)\big)(s-1).
\end{align*}
Since $V^{\Inf_{\mathfrak{t}}^{\mathfrak{g}}(\mu)}=\CC\spanning\{v_{g\Inf_{\mathfrak{t}}^{\mathfrak{g}}(\lambda)}\mid g\in G\}$, this computation implies that $\varphi$ is injective.  

Finally, for $s,t\in T$, $h\in H$, and $\mu=s\Inf_{\mathfrak{t}}^{\mathfrak{g}}(\lambda)$,
$$ht \varphi(v_\mu) = t v_{\Res_{\mathfrak{t}}^{\mathfrak{g}}(\mu)}=\vartheta\big(\mu(t^{-1}-1)\big) v_{t\Res_{\mathfrak{t}}^{\mathfrak{g}}(\mu)}=
\vartheta\big(\mu(t^{-1}-1)\big) v_{\Res_{\mathfrak{t}}^{\mathfrak{g}}(ht\mu)}=\varphi(ht v_\mu),$$
so $\varphi$ is a $G$-module isomorphism.
\end{proof}

We will be primarily be interested in the superinflation function between parabolic subgroups of $U_n(q)$.  In this case, if $U_K(q)\subseteq U_L(q)$,  then 
$$\Sinf_{U_K(q)}^{U_L(q)}(\chi^\lambda)=\chi^{\langle\lambda\rangle_L}.$$
For example,
$$\begin{array}{ccccc} U_{\{2,3,5,7\}} & \overset{\Sinf}{\longrightarrow} & U_{\{1\}\cup\{2,3,5,7\}} & \overset{\Sinf}{\longrightarrow} & U_7\\ 
\chi^{\xymatrix@C=.3cm{*={\scs\circ} & *={\scs\bullet}  \ar @{-} @/^{.25cm}/ [rrr]^{a} &*={\scs\bullet} & *={\scs\circ} & *={\scs\bullet} \ar @{-} @/^{.2cm}/ [rr]^{b} & *={\scs\circ} & *={\scs\bullet}}} & \mapsto &
\chi^{\xymatrix@C=.3cm{*={\scs\bullet} & *={\scs\bullet}  \ar @{-} @/^{.25cm}/ [rrr]^{a} &*={\scs\bullet} &*={\scs\circ} &*={\scs\bullet} \ar @{-} @/^{.2cm}/ [rr]^{b} &*={\scs\circ} &*={\scs\bullet}}}
 & \mapsto & 
 \chi^{\xymatrix@C=.3cm{*={\scs\bullet} & *={\scs\bullet}  \ar @{-} @/^{.25cm}/ [rrr]^{a} &*={\scs\bullet} &*={\scs\bullet} &*={\scs\bullet} \ar @{-} @/^{.2cm}/ [rr]^{b} &*={\scs\bullet} &*={\scs\bullet}}}. \end{array}$$
Thus, superinflation allows us to embed $C_m(q)\subseteq C_n(q)$ for all $m<n$, although this embedding still depends on the embedding of $U_m(q)$ inside $U_n(q)$.

\vspace{,25cm}

\noindent\textbf{Remarks.}  \hfill
\begin{enumerate}
\item While the inflation function does match up with the usual inflation when possible, it does not generally behave as nicely as the usual inflation function.  In particular, it is no longer generally true that $\Res_T^G\circ \Sinf_T^G(\chi)=\chi$ for $\chi$ a class function of $T$.  For example,   
 $$\chi^{\xymatrix@C=.3cm{*={\scs\circ} & *={\scs\bullet}  \ar @{-} @/^{.25cm}/ [rrr]^{a} &*={\scs\bullet} & *={\scs\circ} & *={\scs\bullet} \ar @{-} @/^{.2cm}/ [rr]^{b} & *={\scs\circ} & *={\scs\bullet}}} (1)=q^1\neq q^3= \chi^{\xymatrix@C=.3cm{*={\scs\bullet} & *={\scs\bullet}  \ar @{-} @/^{.25cm}/ [rrr]^{a} &*={\scs\bullet} &*={\scs\bullet} &*={\scs\bullet} \ar @{-} @/^{.2cm}/ [rr]^{b} &*={\scs\bullet} &*={\scs\bullet}}}(1).$$
\item Suppose $K=K_1\cup K_2\in \cS_{m+n}$ with $|K_1|=m$ and $|K_2|=n$.  For $\lambda\in \cS_m(q)$ and $\mu\in \cS_n(q)$, 
$$\Sinf_{U_m\times_K U_n}^{U_{m+n}}(\chi^\lambda\times \chi^\mu)=\chi^{\lambda\cup_K \mu}.$$
\end{enumerate}

%
%

\subsection{Restrictions}\label{SectionRestrictions}

In this section we give algorithms for computing restrictions between parabolic subgroups of $U_n(q)$.   Since supercharacters decompose into tensor products of arcs, for $\lambda\in \cS_n(q)$,
$$\chi^\lambda=\prod_{i\larc{a}l\in A(\lambda)} \chi^{\langle i\larc{a}l\rangle_n},$$
our strategy is to compute restrictions to for each $\chi^{\langle i\larc{a}l\rangle_n}$.  We then use a tensor product result in Section \ref{SectionTensorProduct} to glue back together the resulting restrictions.

We begin with two observations, and then Proposition \ref{FirstIteration} and Theorem \ref{RestrictionRule}  combine to give a general algorithm.  Recall that for $K=K_1\cup K_2\cup\cdots\cup K_\ell\in \cS_n$, $U_K$ is a subgroup of $U_n(q)$ isomorphic to
$$U_{|K_1|}\times U_{|K_2|}\times\cdots \times U_{|K_\ell|}.$$

\begin{proposition}\label{FirstSplitProposition}
Let $U_K\subseteq U_L$ be parabolic subgroups of $U_n$ with $L=L_1\cup L_2\cup\cdots\cup L_\ell\in \cS_n$.  Then
$$\Res^{U_L}_{U_K}(\chi^{\lambda_1}\times\cdots\times \chi^{\lambda_\ell})=\Res^{U_{L_1}}_{U_{K_1}}(\chi^{\lambda_1})\times \Res^{U_{L_2}}_{U_{K_2}}(\chi^{\lambda_2})\times \cdots \times \Res^{U_{L_\ell}}_{U_{K_\ell}}(\chi^{\lambda_\ell})$$
where $U_{K_j}$ is the parabolic subgroup of $U_{L_j}$ corresponding to the vertices $L_j$.
\end{proposition}

The next proposition gives information about each factor in Proposition \ref{FirstSplitProposition}.

\begin{proposition}\label{LeviRestrictionRule}
For $i<l$, $a\in \FF_q^\times$ and $K=K_1\cup K_2\cup\ldots\cup K_\ell\in \cS_n$,
$$\Res_{U_K}^{U_n}(\chi^{\langle i\larc{a} l \rangle_n})= \frac{\Res_{U_{K_1}}^{U_n}(\chi^{\langle i\larc{a} l \rangle_n})}{q^{|\{i<k<l\mid k\notin K_1\}|}}\times \frac{\Res_{U_{K_2}}^{U_n}(\chi^{\langle i\larc{a} l \rangle_n})}{q^{|\{i<k<l\mid k\notin K_2\}|}}\times \cdots \times \frac{\Res_{U_{K_\ell}}^{U_n}(\chi^{\langle i\larc{a} l \rangle_n})}{q^{|\{i<k<l\mid k\notin K_\ell\}|}}$$
\end{proposition}
\begin{proof}
It follows from (\ref{CharacterFormula}) (see also \cite{TV07p} for a more general result) that we can factor the character values across the direct product as
\begin{align*}
\Res_{U_K}^{U_n}(\chi^{\langle i\larc{a} l \rangle_n})(u_{\mu^{(1)}}, \ldots , u_{\mu^{(\ell)}}) 
&=\chi^{\langle i\larc{a} l \rangle_n}(1)\prod_{j=1}^\ell \frac{\chi^{\langle i\larc{a} l \rangle_n}(u_{\mu^{(j)}})}{\chi^{\langle i\larc{a} l \rangle_n}(1)}\\
&=\chi^{\langle i\larc{a} l \rangle_n}(1)\prod_{j=1}^\ell \frac{\Res^{U_n}_{U_{K_j}}(\chi^{\langle i\larc{a} l \rangle_n})(u_{\mu^{(j)}})}{\Res^{U_n}_{U_{K_j}}(\chi^{\langle i\larc{a} l \rangle_n})(1)},\\
\intertext{and by (\ref{DegreeFormula}),}
&=q^{l-i-1}\prod_{j=1}^\ell \frac{\Res^{U_n}_{U_{K_j}}(\chi^{\langle i\larc{a} l \rangle_n})(u_{\mu^{(j)}})}{q^{l-i-1}}\\
&=\frac{1}{q^{(\ell-1)(l-i-1)}}\prod_{j=1}^\ell\Res^{U_n}_{U_{K_j}}(\chi^{\langle i\larc{a} l \rangle_n})(u_{\mu^{(j)}})\\
&= \prod_{j=1}^\ell \frac{\Res_{U_{K_j}}^{U_n}(\chi^{\langle i\larc{a} l \rangle_n})(u_{\mu^{(j)}})}{q^{|\{i<k<l\mid k\notin K_j\}|}},
\end{align*}
as desired.
\end{proof}

To compute restrictions, we first consider very specific subgroups of $U_n(q)$.  For $1<j<k$, let
$$[j,k]=\{j,j+1,\ldots, k-1,k\},$$
and 
$$U_{[i,l]}=\{u\in U_n\mid u_{jk}\neq 0 \text{ implies } i\leq j\leq k\leq l\}.$$

\begin{proposition}\label{FirstIteration} 
Let $1\leq j<k\leq n$ and $S=[j,k]$.  For $1\leq i<l\leq n$ and $a\in \FF_q^\times$,
$$\Res^{U_n}_{U_{[j,k]}}(\chi^{\langle i \larc{a} l \rangle_n})=\left\{\begin{array}{ll} \chi^{\langle i \larc{a} l \rangle_S}, & \text{if $j\leq i<l\leq k$,}\\
\dd q^{j-i-1}\bigg(\One+\sum_{j\leq j'<l\atop b\in \FF_q^\times} \chi^{\langle j' \larc{b} l \rangle_S} \bigg), & \text{if $i<j<l\leq k$}\\
\dd q^{l-k-1}\bigg(\One+\sum_{i\leq k'\leq k\atop b\in \FF_q^\times} \chi^{\langle i\larc{b} k' \rangle_S} \bigg), & \text{if $j\leq i <k<l$}\\
\dd \frac{q^{l-i-1}}{q^{k-j+1}}\bigg((k-j)(q-1)+q)\One+(q-1)\hspace{-.25cm}\sum_{j\leq j'<k'\leq k\atop b\in \FF_q^\times} \chi^{\langle j'\larc{b} k' \rangle_S}\bigg), & \text{if $i<j<k<l$,}\\
q^{l-i-1}\One, & \text{otherwise}.\end{array} 
\right. $$
\end{proposition}
\begin{proof}
Cases 1 and 5 follow directly from  (\ref{CharacterFormula}). 

For Cases 2, first assume that $i=1$, $j=2$ and $l\leq k=n$.  Then \cite{TV07} implies the result ($q^{j-i-1}=1$ in this special case). 

To obtain the remaining cases we iterate by removing one column or row at a time.  For $i<j< l\leq k\leq  n$,
$$\Res^{U_n}_{U_{[j,k]}}(\chi^{\langle i \larc{a} l \rangle_n})= \Res^{U_{[j,n]}}_{U_{[j,k]}}\Res_{U_{[j,n]}}^{U_n}(\chi^{\langle i \larc{a} l \rangle_n}),$$
where
\begin{align*}
\Res_{U_{[j,n]}}^{U_{n}} &(\chi^{\langle i\larc{a}l \rangle_n}) =\Res_{U_{[j,n]}}^{U_{[j-1,n]}}\circ\cdots\circ \Res_{U_{[i+2,n]}}^{U_{[i+1,n]}}\circ \Res_{U_{[i+1,n]}}^{U_n}(\chi^{\langle i\larc{a}l\rangle_n})\\
&=\Res_{U_{[j,n]}}^{U_{[j-1,n]}}\circ\cdots\circ \Res_{U_{[i+2,n]}}^{U_{[i+1,n]}}\bigg(\One+\sum_{i+1\leq j'<l\atop b\in \FF_q^\times}\chi^{\langle j'\larc{b}l \rangle_{[i+1,n]}}\bigg)\\
&=\Res_{U_{[j,n]}}^{U_{[j-1,n]}}\circ\cdots\circ \Res_{U_{[i+3,n]}}^{U_{[i+2,n]}}\bigg(\One+\hspace{-.4cm}\sum_{i+2\leq j'<l\atop b\in \FF_q^\times}\chi^{\langle j'\larc{b}l\rangle_{[i+1,n]}}+\sum_{b\in \FF_q^\times}\Res_{U_{[i+2,n]}}^{U_{[i+1,n]}}(\chi^{\langle i+1\larc{b}n \rangle_{[i+1,n]}})\bigg)\\
&=\Res_{U_{[j,n]}}^{U_{[j-1,n]}}\circ\cdots\circ \Res_{U_{[i+3,n]}}^{U_{[i+2,n]}}\bigg(q\One+q\sum_{i+2\leq j'<l\atop b\in \FF_q^\times}\chi^{\langle j'\larc{b}l\rangle_{[i+1,n]}}\bigg)\\
&=q\Res_{U_{[j,n]}}^{U_{[j-1,n]}}\circ\cdots\circ \Res_{U_{[i+3,n]}}^{U_{[i+2,n]}}\bigg(\One+\sum_{i+2\leq j'<l\atop b\in \FF_q^\times}\chi^{\langle j'\larc{b}l\rangle_{[i+1,n]}}\bigg)\\
\intertext{and iterate to obtain}
&=q^{j-i-1}\One+q^{j-i-1}\sum_{j\leq j'<l\atop b\in \FF_q^\times}\chi^{\langle j'\larc{b}l\rangle_S},
\end{align*}
giving Case 2.  Case 3 follows by a symmetric argument.

If $i<j<k<l$, then by Case 2,
\begin{align*}
&\Res_{U_{[j,k]}}^{U_{n}} (\chi^{\langle i\larc{a}l\rangle_n}) = \Res_{U_{[j,k]}}^{U_{[j,n]}}\circ \Res_{U_{[j,n]}}^{U_{n}}(\chi^{\langle i\larc{a}l\rangle_n})\\
&= \Res_{U_{[j,k]}}^{U_{[j,n]}} \bigg(q^{j-i-1}\One+q^{j-i-1}\sum_{j\leq j'< l\atop b\in \FF_q^\times}\chi^{\langle j'\larc{b}l\rangle_{[j,n]}}\bigg)\\
&= \Res_{U_{[j,k]}}^{U_{[j,n]}} \bigg(q^{j-i-1}\One+q^{j-i-1}\sum_{j\leq j'< k\atop b\in \FF_q^\times}\chi^{\langle j'\larc{b}l\rangle_{[j,n]}} +q^{j-i-1}\sum_{k\leq j'< l\atop b\in \FF_q^\times}\chi^{\langle j'\larc{b}l\rangle_{[j,n]}}\bigg),\\
\intertext{and by Case 3,}
&=q^{j-i-1}\One+q^{j-i-1}\sum_{j\leq j'<k\atop b\in \FF_q^\times}\bigg(q^{l-k-1}\One+q^{l-k-1}\sum_{j'<k'\leq k\atop c\in \FF_q^\times}\chi^{\langle j'\larc{c}k'\rangle_S}\bigg)+q^{j-i-1}(q-1)\sum_{k\leq j'<l} q^{l-j'-1}\One\\
&=q^{j-i-1}\One+q^{l-k+j-i-2}(k-j)(q-1)\One+q^{l-k+j-i-2}\sum_{j\leq j'<k'\leq k\atop b,c\in \FF_q^\times}\chi^{\langle j'\larc{c}k'\rangle_S}+q^{j-i-1}(q^{l-k}-1)\One\\
&=(q^{l-k+j-i-2}(q-1)(k-j)+q^{l-k+j-i-1})\One +q^{l-k+j-i-2}(q-1)\sum_{j\leq j'<k'\leq k\atop c\in \FF_q^\times}\chi^{\langle j'\larc{c}k'\rangle_S}\\
&=\frac{q^{l-i-1}}{q^{k-j+1}} \bigg(((k-j)(q-1)+q)\One+(q-1)\sum_{j\leq j'<k'\leq k\atop c\in \FF_q^\times}\chi^{\langle j'\larc{c}k'\rangle_S}\bigg),
\end{align*}
giving Case 4.
\end{proof}

\noindent\textbf{Example.} Let $n=7$, $j=2$, $k=5$, so that
$$U_{[2,5]}=\left\{\left(\begin{smallmatrix} 
1 & 0 & 0 & 0 & 0 & 0 & 0\\ 
0 & 1 & \ast & \ast & \ast & 0 & 0 \\
0 & 0 & 1 & \ast & \ast  & 0 & 0\\
0 & 0 & 0 & 1 & \ast  & 0 & 0\\
0 & 0 & 0 & 0 & 1 & 0 & 0\\
0 & 0 & 0 & 0 & 0 & 1 & 0\\ 
0 & 0 & 0 & 0 & 0 & 0 & 1\end{smallmatrix}\right)\right\} \subseteq
\left\{\left(\begin{smallmatrix} 
1 & \ast & \ast & \ast & \ast & \ast & \ast\\ 
0 & 1 & \ast & \ast & \ast & \ast & \ast \\
0 & 0 & 1 & \ast & \ast  & \ast & \ast\\
0 & 0 & 0 & 1 & \ast  & \ast & \ast\\
0 & 0 & 0 & 0 & 1 & \ast & \ast\\
0 & 0 & 0 & 0 & 0 & 1 & \ast\\ 
0 & 0 & 0 & 0 & 0 & 0 & 1\end{smallmatrix}\right)\right\}=U_7$$
Then 
\begin{align*}
\Res_{U_{[2,5]}}^{U_7}(\chi^{\xymatrix@C=.25cm{*={\scs\bullet} & *={\scs\bullet}  \ar @{-} @/^{.25cm}/ [rrr]^{a} &*={\scs\bullet} &*={\scs\bullet} &*={\scs\bullet} &*={\scs\bullet} &*={\scs\bullet} }})&=\chi^{\xymatrix@C=.25cm{*={\scs\circ} & *={\scs\bullet}  \ar @{-} @/^{.25cm}/ [rrr]^{a} &*={\scs\bullet} &*={\scs\bullet} &*={\scs\bullet} &*={\scs\circ} &*={\scs\circ}}}\\
\Res_{U_{[2,5]}}^{U_7}(\chi^{\xymatrix@C=.25cm{*={\scs\bullet} \ar @{-} @/^{.25cm}/ [rrrr]^{a} & *={\scs\bullet} &*={\scs\bullet} &*={\scs\bullet} &*={\scs\bullet} &*={\scs\bullet} &*={\scs\bullet} }})&=
\chi^{\xymatrix@C=.25cm{*={\scs\circ} & *={\scs\bullet} &*={\scs\bullet} &*={\scs\bullet} &*={\scs\bullet} &*={\scs\circ} &*={\scs\circ} }} + \sum_{b\in \FF_q^\times}\chi^{\xymatrix@C=.25cm{*={\scs\circ} & *={\scs\bullet}  \ar @{-} @/^{.25cm}/ [rrr]^{b}&*={\scs\bullet} &*={\scs\bullet} &*={\scs\bullet} &*={\scs\circ} &*={\scs\circ} }}
+ \sum_{b\in \FF_q^\times}\chi^{\xymatrix@C=.25cm{*={\scs\circ} & *={\scs\bullet}&*={\scs\bullet}   \ar @{-} @/^{.25cm}/ [rr]^{b}&*={\scs\bullet} &*={\scs\bullet} &*={\scs\circ} &*={\scs\circ} }}
+ \sum_{b\in \FF_q^\times}\chi^{\xymatrix@C=.25cm{*={\scs\circ} & *={\scs\bullet}&*={\scs\bullet}  &*={\scs\bullet} \ar @{-} @/^{.25cm}/ [r]^{b} &*={\scs\bullet} &*={\scs\circ} &*={\scs\circ} }}\\
\Res_{U_{[2,5]}}^{U_7}(\chi^{\xymatrix@C=.25cm{*={\scs\bullet} \ar @{-} @/^{.25cm}/ [rrrrrr]^{a} & *={\scs\bullet} &*={\scs\bullet} &*={\scs\bullet} &*={\scs\bullet} &*={\scs\bullet} &*={\scs\bullet} }})
&=
q\bigg((4q-3)\chi^{\xymatrix@C=.25cm{*={\scs\circ} & *={\scs\bullet} &*={\scs\bullet} &*={\scs\bullet} &*={\scs\bullet} &*={\scs\circ} &*={\scs\circ} }}
+(q-1)\sum_{b\in \FF_q^\times}\chi^{\xymatrix@C=.25cm{*={\scs\circ} & *={\scs\bullet}  \ar @{-} @/^{.25cm}/ [r]^{b}&*={\scs\bullet} &*={\scs\bullet} &*={\scs\bullet} &*={\scs\circ} &*={\scs\circ} }}
+(q-1) \sum_{b\in \FF_q^\times}\chi^{\xymatrix@C=.25cm{*={\scs\circ} & *={\scs\bullet}&*={\scs\bullet}   \ar @{-} @/^{.25cm}/ [rr]^{b}&*={\scs\bullet} &*={\scs\bullet} &*={\scs\circ} &*={\scs\circ} }}\\
&+ (q-1)\sum_{b\in \FF_q^\times}\chi^{\xymatrix@C=.25cm{*={\scs\circ} & *={\scs\bullet}&*={\scs\bullet}  &*={\scs\bullet} \ar @{-} @/^{.25cm}/ [r]^{b} &*={\scs\bullet} &*={\scs\circ} &*={\scs\circ} }}+(q-1)\sum_{b\in \FF_q^\times}\chi^{\xymatrix@C=.25cm{*={\scs\circ} & *={\scs\bullet} \ar @{-} @/^{.25cm}/ [rrr]^{b}&*={\scs\bullet}  &*={\scs\bullet}  &*={\scs\bullet} &*={\scs\circ} &*={\scs\circ} }}+(q-1)\sum_{b\in \FF_q^\times}\chi^{\xymatrix@C=.25cm{*={\scs\circ} & *={\scs\bullet} \ar @{-} @/^{.25cm}/ [rr]^{b}&*={\scs\bullet}  &*={\scs\bullet}  &*={\scs\bullet} &*={\scs\circ} &*={\scs\circ} }}\\
&+(q-1)\sum_{b\in \FF_q^\times}\chi^{\xymatrix@C=.25cm{*={\scs\circ} & *={\scs\bullet} &*={\scs\bullet} \ar @{-} @/^{.25cm}/ [r]^{b} &*={\scs\bullet}  &*={\scs\bullet} &*={\scs\circ} &*={\scs\circ} }}\bigg)\\
\Res_{U_{[2,5]}}^{U_7}(\chi^{\xymatrix@C=.25cm{*={\scs\bullet} & *={\scs\bullet} &*={\scs\bullet} &*={\scs\bullet} &*={\scs\bullet}  \ar @{-} @/^{.25cm}/ [rr]^{a}  &*={\scs\bullet} &*={\scs\bullet} }})&=q\chi^{\xymatrix@C=.25cm{*={\scs\circ} & *={\scs\bullet} &*={\scs\bullet} &*={\scs\bullet} &*={\scs\bullet} &*={\scs\circ} &*={\scs\circ} }}
\end{align*}

Theorem \ref{RestrictionRule}, below, generalizes Proposition \ref{FirstIteration} to the required level of generality.  For $S\subseteq \{1,2,\ldots, n\}$, let
$$U_S=\{u\in U_n\mid u_{ij}\neq 0 \text{ implies } i,j\in S\}.$$
Note that while $U_S$ is not itself a parabolic subgroup of $U_n$, it is isomorphic to the parabolic subgroup $U_{\langle S\rangle_n}$.

\begin{thm}\label{RestrictionRule}
Let $S\subseteq \{1,2,\ldots, n\}$.  Then for $1\leq i<l\leq n$ and $a\in \FF_q^\times$,
$$\Res^{U_n}_{U_S}(\chi^{\langle i\larc{a}l\rangle_n})=\left\{\begin{array}{ll}
\dd q^{|\{i<k<l\mid k\notin S\}|} \chi^{\langle i\larc{a}l\rangle_S}, & \text{if $i,l\in S$,}\\
\dd q^{|\{i<k<l\mid k\notin S\}|}\bigg(\One+\sum_{i<j<l, j\in S\atop b\in \FF_q^\times} \chi^{\langle j\larc{b} l \rangle_S}\bigg),& \text{if $i\notin S, l\in S$,}\\
\dd q^{|\{i<k<l\mid k\notin S\}|}\bigg(\One+\sum_{i<k<l, k\in S\atop b\in \FF_q^\times} \chi^{\langle i\larc{b} k \rangle_S}\bigg),& \text{if $i\in S, l\notin S$,}\\
\dd q^{|\{i<k'<l \mid k'\notin S\}|}\bigg((|S\cap[i,l]|(q-1)+1)\One+(q-1)\hspace{-.5cm}\sum_{i<j'< k'<l\atop j',k'\in S, c\in \FF_q^\times} \hspace{-.3cm}\chi^{\langle j'\larc{c} k' \rangle_S}\bigg), &  \text{if $i,l\notin S$.}
\end{array}\right.$$
\end{thm}
\begin{proof}
Note that it follows from the character formulas (\ref{CharacterFormula}) and (\ref{DegreeFormula}) that if $i,l\in S$, then
$$\frac{\Res^{U_n}_{U_S}(\chi^{\langle i\larc{a} l \rangle_n})}{\chi^{\langle i\larc{a} l \rangle_n}(1)}=\frac{\chi^{\langle i\larc{a} l \rangle_S}}{\chi^{\langle i\larc{a} l \rangle_S}(1)}.$$
Therefore,
\begin{equation}\label{Case1Equation}
\Res^{U_n}_{U_S}(\chi^{\langle i\larc{a} l \rangle_n})=\frac{\chi^{\langle i\larc{a} l \rangle_n}(1)}{\chi^{\langle i\larc{a} l \rangle_S}(1)}\chi^{\langle i\larc{a} l \rangle_S}=q^{|\{i<k<l\mid k\notin S\}|}\chi^{\langle i\larc{a} l \rangle_S}
\end{equation}
gives the first case.

For Case 2, assume that $i\notin S$, and $l\in S$.  By Proposition \ref{FirstIteration}, Case 2 and then (\ref{Case1Equation}),
\begin{align}
\Res^{U_n}_{U_S}(\chi^{\langle i\larc{a} l \rangle_n})&=\One+\sum_{i<j<l\atop b\in \FF_q^\times} \Res_{U_S}^{U_{[i+1,l]}}(\chi^{\langle j\larc{b} l \rangle_{[i+1,l]}})\notag\\
&= \One + \sum_{i<j<l, j\in S \atop b\in \FF_q^\times} q^{|\{j<k<l \mid k\notin S\}|}\chi^{\langle j\larc{b} l \rangle_S} + \sum_{i<j<l,j\notin S\atop b\in \FF_q^\times}  \Res_{U_S}^{U_{[i+1,l]}}(\chi^{\langle j\larc{b} l \rangle_{[i+1,l]}}).\label{RestrictionComputationOne}
\end{align}
If $j'$ is minimal such that $i<j'<l$ and $j'\notin S$, then by Proposition \ref{FirstIteration}, Case 2,
\begin{align*}
\sum_{i<j<l,j\notin S\atop b\in \FF_q^\times} & \Res_{U_S}^{U_{[i+1,l]}}(\chi^{\langle j\larc{b} l \rangle_{[i+1,l]}})\\
&= \sum_{b\in \FF_q^\times}\Res_{U_S}^{U_{[j'+1,l]}}\Res_{U_{[j'+1,l]}}^{U_{[i+1,l]}}(\chi^{\langle j'\larc{b}l\rangle_{[i+1,l]}})+\sum_{j'<j<l,j\notin S\atop b\in \FF_q^\times}  \Res_{U_S}^{U_{[j'+1,l]}}(\chi^{\langle j\larc{b} l \rangle_{[j'+1,l]}})\\
&= (q-1)\One + (q-1) \sum_{j'<j<l\atop c\in \FF_q^\times} \Res_{U_S}^{U_{[j'+1,l]}}(\chi^{\langle j \larc{c} l \rangle_{[j'+1,l]}})+\sum_{j'<j<l,j\notin S\atop b\in \FF_q^\times} \Res_{U_S}^{U_{[j'+1,l]}}(\chi^{\langle j\larc{b} l \rangle_{[j'+1,l]}})\\
&=(q-1)\One + (q-1)\sum_{j'<j<l,j\in S\atop c\in \FF_q^\times} q^{|\{j<k<l \mid k\notin S\}|} \chi^{\langle j \larc{c} l \rangle_S} + q\sum_{j'<j<l,j\notin S\atop c\in \FF_q^\times} \Res_{U_S}^{U_{[j'+1,l]}}(\chi^{\langle j \larc{c} l \rangle_{[j'+1,l]}}).
\end{align*}
Applying this equality to (\ref{RestrictionComputationOne}),
\begin{align*}
\Res^{U_n}_{U_S}(\chi^{\langle i\larc{a} l \rangle_n})=q\One &+q^{|\{i<k<l \mid k\notin S\}|}\sum_{i<j<j'\atop b\in \FF_q^\times} \chi^{\langle j\larc{b} l \rangle_S} + q\sum_{j'<j<l, j\in S\atop b\in \FF_q^\times}  q^{|\{j<k<l \mid k\notin S\}|} \chi^{\langle j\larc{b} l \rangle_S} \\
&+ q\sum_{j'<j<l,j\notin S\atop c\in \FF_q^\times} \Res_{U_S}^{U_{[j'+1,l]}}(\chi^{\langle j \larc{c} l \rangle_{[j'+1,l]}})\\
=q\One&+q^{|\{i<k<l \mid k\notin S\}|}\sum_{i<j<j''\atop b\in \FF_q^\times} \chi^{\langle j\larc{b} l \rangle_S} + q\bigg(\sum_{j''<j<l, j\in S\atop b\in \FF_q^\times}  q^{|\{j<k<l \mid k\notin S\}|} \chi^{\langle j\larc{b} l \rangle_S} \\
&+ \sum_{j'<j<l,j\notin S\atop c\in \FF_q^\times} \Res_{U_S}^{U_{[j'+1,l]}}(\chi^{\langle j \larc{c} l \rangle_{[j'+1,l]}})\bigg).
\end{align*}
where $j''$ is minimal such that $j'<j''<l$ and $j''\notin S$.  Iterating this process gives Case 2.  A symmetric argument also gives Case 3.

Suppose $i,l\notin S$ with $|S\cap[i,l]|\neq 0$.  Let $k'<l$ maximal such that $k'\in S$.  We apply Cases 2 and 3 consecutively, to deduce Case 3. By Case 2,
\begin{align*}
&\Res^{U_n}_{U_S}(\chi^{\langle i\larc{a} l \rangle_n})=\Res_{U_S}^{U_{S\cup\{l\}}}\circ \Res_{U_{S\cup\{l\}}}^{U_n}(\chi^{\langle i\larc{a} l \rangle_n})\\
&=q^{|\{i<k<l \mid k\notin S\}|}\Res_{U_S}^{U_{S\cup\{l\}}}\bigg(\One+\sum_{i<j<l, j\in S\atop b\in \FF_q^\times} \chi^{\langle j\larc{b} l \rangle_{S\cup\{l\}}}\bigg)\\
&=q^{|\{i<k<l \mid k\notin S\}|}\bigg(\One+\sum_{i<j<l, j\in S\atop b\in \FF_q^\times} \Res_{U_S}^{U_{S\cup\{l\}}}(\chi^{\langle j \larc{b} l \rangle_{S\cup\{l\}}})\bigg)\\
&=q^{|\{i<k<l \mid k\notin S\}|}\bigg(\One+\sum_{i<j<k', j\in S\atop b\in \FF_q^\times} \Res_{U_S}^{U_{S\cup\{l\}}}(\chi^{\langle j \larc{b} l \rangle_{S\cup\{l\}}})+\sum_{b\in \FF_q^\times} \Res_{U_S}^{U_{S\cup\{l\}}}(\chi^{\langle k' \larc{b} l \rangle_{S\cup\{l\}}})\bigg)\\
\intertext{and then  by Case 3,}
&=q^{|\{i<k<l \mid k\notin S\}|}\bigg(\One+|S\cap[i,k')|(q-1)\One+(q-1)\sum_{i<j< k<l, j,k\in S\atop c\in \FF_q^\times}\chi^{\langle j\larc{c} k\rangle_S}+(q-1)\One\bigg)\\
&=q^{|\{i<k<l \mid k\notin S\}|}\bigg((|S\cap[i,k)|(q-1)+q)\One+(q-1)\sum_{i<j< k<l, j,k\in S\atop c\in \FF_q^\times} \chi^{\langle j\larc{c} k \rangle_S}\bigg)\\
&=q^{|\{i<k<l \mid k\notin S\}|}\bigg(((|S\cap[i,l]|-1)(q-1)+q)\One+(q-1)\sum_{i<j< k<l, j,k\in S\atop c\in \FF_q^\times} \chi^{\langle j\larc{c} k \rangle_S}\bigg)\\
&=q^{|\{i<k<l \mid k\notin S\}|}\bigg((|S\cap[i,l]|(q-1)+1)\One+(q-1)\sum_{i<j< k<l, j,k\in S\atop c\in \FF_q^\times} \chi^{\langle j\larc{c} k \rangle_S}\bigg).
\end{align*}
On the other hand, if $i,l\notin S$ and $|S\cap[i,l]|= 0$, then there exists $j<k$ such that $S\subseteq [j,k]$ and $[i,l]\cap[j,k]=\emptyset$.  Thus, by Proposition \ref{FirstIteration}, Case 5,
\begin{align*}
\Res^{U_n}_{U_S}(\chi^{\langle i\larc{a} l \rangle_n})&=\Res_{U_S}^{U_{[j,k]}}\circ \Res_{U_{[j,k]}}^{U_n}(\chi^{\langle i\larc{a} l \rangle_n})\\
&=q^{l-i-1}\One\\
&=q^{|\{i<k<l \mid k\notin S\}|}\bigg((|S\cap[i,l]|(q-1)+1)\One+(q-1)\sum_{i<j< k<l, j,k\in S\atop c\in \FF_q^\times} \chi^{\langle j\larc{c} k \rangle_S}\bigg),
\end{align*}
where the sum is empty in this case.
\end{proof}

%
%

\subsection{Tensor products} \label{SectionTensorProduct}

We have seen in the previous section that when we decompose supercharacters into tensor products of irreducible characters, the restriction rules are manageable to compute.  This section explains how to glue back together the resulting products of characters.  The main result -- Corollary \ref{TensorResult} -- has been computed by Andr\'e for large primes in \cite[Lemmas 6--8]{An95} and for arbitrary primes by Yan in \cite[Propositions 7.2-7.5]{Ya01}, but we reprove it here quickly using the machinery developed in this paper.

We begin with a lemma that further establishes the relationship between tensor products and restrictions.  For $H\subseteq G$ and $\chi$ a superclass function of $G$, let 
$$\mathrm{SinfRes}_{H}^G(\chi)=\Sinf_{H}^G\Res_{H}^G(\chi).$$

\begin{lemma}\label{TensorLemma} For $i<j<k<l$, 
\begin{align*}
\chi^{\langle i \larc{a} k \rangle_n}\otimes\chi^{\langle i \larc{b} l \rangle_n} &= \mathrm{SinfRes}_{U_{[i+1,l]}}^{U_{n}}(\chi^{\langle i \larc{a} k \rangle_n})\otimes\chi^{\langle i \larc{b} l \rangle_n}, & & a,b\in \FF_q^\times,\\
\chi^{\langle i\larc{a} l \rangle_n}\otimes\chi^{\langle j\larc{b} l \rangle_n} &= \chi^{\langle i\larc{a} l \rangle_n}\otimes \mathrm{SinfRes}_{U_{[i,l-1]}}^{U_{n}}(\chi^{\langle j\larc{b} l \rangle_n}),& & a,b\in \FF_q^\times,\\ 
\chi^{\langle i\larc{a} l \rangle_n}\otimes\chi^{\langle i \larc{-a} l \rangle_n}&= \mathrm{SinfRes}_{U_{[i+1,l]}}^{U_n}(\chi^{\langle i\larc{a} l \rangle_n})\otimes \mathrm{SinfRes}_{U_{[i,l-1]}}^{U_n}(\chi^{\langle i \larc{-a} l \rangle_n}),& & a\in \FF_q^\times,\\
\chi^{\langle i\larc{a} l \rangle_n}\otimes\chi^{\langle i \larc{b} l \rangle_n} &= \chi^{\langle i \larc{a+b} l \rangle_n}\otimes  \mathrm{SinfRes}_{U_{[i+1,l-1]}}^{U_{n}}(\chi^{\langle i \larc{a+b} l \rangle_n}),& & a,b\in \FF_q^\times, b\neq -a.
\end{align*}
\end{lemma}
\begin{proof}
For the first case, note that by (\ref{CharacterFormula}),
$$\chi^{\langle i \larc{a} k \rangle_n}\otimes\chi^{\langle i \larc{b} l \rangle_n}(u_\mu)=\left\{\begin{array}{ll}\dd\frac{\chi^{\langle i \larc{a} k \rangle_n}(1)\chi^{\langle i \larc{b} l \rangle_n}(1)\vartheta(b\tau_\mu(i\larc{}l))}{q^{|\{j'\larc{}k'\in \mu\mid i<j'<k'<k\}|+|\{j'\larc{}k'\in \mu\mid i<j'<k'<l\}|}} , & \begin{array}{@{}l} \text{if $i<j'<k$ implies}\\ i\larc{}j',j'\larc{}k\notin A(\mu),\\  \text{$i<j'<l$ implies}\\ i\larc{}j',j'\larc{}l\notin A(\mu);\end{array}\\
0, & \text{otherwise}.\end{array}\right.$$
where the $\vartheta(a\tau_\mu(i\larc{}k))$ term is missing since if $i\larc{}k\in A(\mu)$ then the character value is zero.   Now the first case follows from the observation that for $\nu\in \cS_{[i+1,l]}(q)$,
$$\Res_{U_{[i+1,l]}}^{U_{n}}(\chi^{\langle i \larc{a} k \rangle_n})(u_\nu)=\left\{\begin{array}{ll}\dd\frac{\chi^{\langle i \larc{a} k \rangle_n}(1)}{q^{|\{j'\larc{}k'\in \nu\mid i<j'<k'<k\}|}}, & \text{if $i<j'<k$ implies $j'\larc{}k\notin A(\nu)$,}\\ 0, & \text{otherwise.}\end{array}\right.$$
Case 2 follows from a symmetric argument, and Cases 3 and 4 are proved by a similar argument.
\end{proof}

Combining Lemma \ref{TensorResult} with Proposition \ref{FirstIteration} we obtain the following corollary.

\begin{cor}\label{TensorResult}  For $i<k$, $j<l$, $a,b\in \FF_q^\times$, and $\{i,k\}\neq \{j,l\}$,
$$\chi^{\langle i \larc{a} k \rangle_n}\otimes\chi^{\langle j\larc{b} l \rangle_n}=\left\{\begin{array}{ll}
 \chi^{\langle\{i\larc{a}k\}\cup\{j\larc{b}l\}\rangle_n}, & \text{if $\{i,k\}\cap\{j,l\}=\emptyset$,}\\
\chi^{\langle i\larc{a}j\larc{b}l\rangle_n}, & \text{if $i<j=k<l$,}\\
\dd \chi^{\langle i \larc{b} l \rangle_n} +\sum_{i<j'<k\atop c\in \FF_q^\times} \chi^{\langle\{j'\larc{c}k\}\cup\{i\larc{b}l\}\rangle_n}, & \text{if $i=j<k<l$,}\\
\dd \chi^{\langle i\larc{a} l \rangle_n} +\sum_{j<k'<l\atop c\in \FF_q^\times} \chi^{\langle\{i\larc{a}l\}\cup\{j\larc{c}k'\}\rangle_n}, & \text{if $i<j<k=l$,}
\end{array}\right.$$
For $i<l$, $a,b\in \FF_q^\times$,
$$\chi^{\langle i\larc{a} l \rangle_n}\otimes\chi^{\langle i \larc{b} l \rangle_n}=\left\{\begin{array}{ll}
\dd \One+\hspace{-.15cm}\sum_{i< j'< l\atop c\in \FF_q^\times} \chi^{\langle i\larc{c}j'\rangle_n}+\hspace{-.15cm}\sum_{i< k'<l\atop c\in \FF_q^\times} \chi^{\langle k'\larc{c}l\rangle_n} + \hspace{-.15cm}\sum_{i<j',k'<l\atop c,d\in \FF_q^\times} \chi^{\langle\{i\larc{c}j'\}\cup\{k'\larc{d}l\}\rangle_n}, & \text{if $b=-a$,}\\
\dd ((l-i-1)(q-1)+1)\chi^{\langle i \larc{a+b} l \rangle_n} + (q-1)\hspace{-.35cm}\sum_{i<j'<k'<l\atop c\in \FF_q^\times}\hspace{-.25cm} \chi^{\langle\{j'\larc{c}k'\}\cup\{i\larc{a+b}l\}\rangle_n}, & \text{otherwise.} \end{array}\right.$$
\end{cor}

\noindent \textbf{Examples.}  Combinatorially, a product of characters is the superimposition of two set-partition diagrams, such as
$$\chi^{\ 
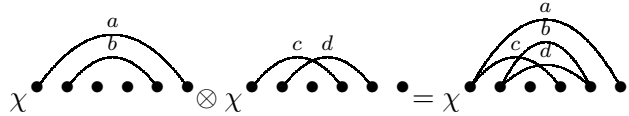
\xymatrix@R=.2cm@C=.4cm{*={\bullet} \ar @{-} @/^{.7cm}/ [rrrrr]^{a} & *={\bullet}  \ar @{-} @/^{.4cm}/ [rrr]^{b}  & *={\bullet} & *={\bullet}  & *={\bullet} & *={\bullet}}}\otimes \chi^{\ \xymatrix@R=.2cm@C=.4cm{*={\bullet}  \ar @{-} @/^{.4cm}/ [rrr]^{c}  & *={\bullet}  \ar @{-} @/^{.4cm}/ [rrr]^{d}  & *={\bullet} & *={\bullet}  & *={\bullet} & *={\bullet}}} =\chi^{\ \xymatrix@R=.2cm@C=.4cm{*={\bullet}  \ar @{-} @/^{.9cm}/ [rrrrr]^{a} \ar @{-} @/^{.4cm}/ [rrr]^{c}  & *={\bullet}  \ar @{-} @/^{.3cm}/ [rrr]^{d}  \ar @{-} @/^{.6cm}/ [rrr]^{b} & *={\bullet} & *={\bullet}  & *={\bullet} & *={\bullet}}} $$  
The tensor product rules then describe how to ``straighten" the resulting diagram.   Usually, we superimpose selectively.  For example,
since
\begin{align*}
\chi^{\ \xymatrix@R=.2cm@C=.4cm{*={\bullet} \ar @{-} @/^{.7cm}/ [rrrrr]^{a} & *={\bullet}  \ar @{-} @/^{.4cm}/ [rrr]^{b}  & *={\bullet} & *={\bullet}  & *={\bullet} & *={\bullet}}}
&=
\chi^{\ \xymatrix@R=.2cm@C=.4cm{*={\bullet} \ar @{-} @/^{.7cm}/ [rrrrr]^{a} & *={\bullet}  & *={\bullet} & *={\bullet}  & *={\bullet} & *={\bullet}}}\otimes \chi^{\ \xymatrix@R=.2cm@C=.4cm{*={\bullet} & *={\bullet}  \ar @{-} @/^{.4cm}/ [rrr]^{b}  & *={\bullet} & *={\bullet}  & *={\bullet} & *={\bullet}}}
\\ 
\chi^{\ \xymatrix@R=.2cm@C=.4cm{*={\bullet}  \ar @{=} @/^{.4cm}/ [rrr]^{c}  & *={\bullet}  \ar @{=} @/^{.4cm}/ [rrr]^{d}  & *={\bullet} & *={\bullet}  & *={\bullet} & *={\bullet}}} 
&=
\chi^{\ \xymatrix@R=.2cm@C=.4cm{*={\bullet}  \ar @{=} @/^{.4cm}/ [rrr]^{c}  & *={\bullet}  & *={\bullet} & *={\bullet}  & *={\bullet} & *={\bullet}}} \otimes \chi^{\ \xymatrix@R=.2cm@C=.4cm{*={\bullet}    & *={\bullet}  \ar @{=} @/^{.4cm}/ [rrr]^{d}  & *={\bullet} & *={\bullet}  & *={\bullet} & *={\bullet}}} 
\end{align*}
 for $d\neq -b$,
\begin{align*}
&\chi^{\ \xymatrix@R=.2cm@C=.4cm{*={\bullet} \ar @{-} @/^{.7cm}/ [rrrrr]^{a} & *={\bullet}  \ar @{-} @/^{.4cm}/ [rrr]^{b}  & *={\bullet} & *={\bullet}  & *={\bullet} & *={\bullet}}}\otimes \chi^{\ \xymatrix@R=.2cm@C=.4cm{*={\bullet}  \ar @{=} @/^{.4cm}/ [rrr]^{c}  & *={\bullet}  \ar @{=} @/^{.4cm}/ [rrr]^{d}  & *={\bullet} & *={\bullet}  & *={\bullet} & *={\bullet}}}\\
&= 
\chi^{\ \xymatrix@R=.2cm@C=.4cm{*={\bullet} \ar @{-} @/^{.7cm}/ [rrrrr]^{a} & *={\bullet}  & *={\bullet} & *={\bullet}  & *={\bullet} & *={\bullet}}}\otimes \chi^{\ \xymatrix@R=.2cm@C=.4cm{*={\bullet} & *={\bullet}  \ar @{-} @/^{.4cm}/ [rrr]^{b}  & *={\bullet} & *={\bullet}  & *={\bullet} & *={\bullet}}}\otimes  \chi^{\ \xymatrix@R=.2cm@C=.4cm{*={\bullet}  \ar @{=} @/^{.4cm}/ [rrr]^{c}  & *={\bullet}  & *={\bullet} & *={\bullet}  & *={\bullet} & *={\bullet}}} \otimes \chi^{\ \xymatrix@R=.2cm@C=.4cm{*={\bullet}    & *={\bullet}  \ar @{=} @/^{.4cm}/ [rrr]^{d}  & *={\bullet} & *={\bullet}  & *={\bullet} & *={\bullet}}} \\
&= 
\chi^{\ \xymatrix@R=.2cm@C=.4cm{*={\bullet} \ar @{-} @/^{.7cm}/ [rrrrr]^{a} \ar @{=} @/^{.4cm}/ [rrr]^(.7){c}  & *={\bullet}  & *={\bullet} & *={\bullet}  & *={\bullet} & *={\bullet}}}\otimes \chi^{\ \xymatrix@R=.2cm@C=.4cm{*={\bullet} & *={\bullet}  \ar @{-} @/^{.3cm}/ [rrr]^{b}  \ar @{=} @/^{.7cm}/ [rrr]^{d} & *={\bullet} & *={\bullet}  & *={\bullet} & *={\bullet}}}\\
&=
 \bigg(\chi^{\ \xymatrix@R=.2cm@C=.4cm{*={\bullet} \ar @{-} @/^{.7cm}/ [rrrrr]^{a} & *={\bullet}  & *={\bullet} & *={\bullet}  & *={\bullet} & *={\bullet}}}+\sum_{e\in \FF_q^\times}\chi^{\ \xymatrix@R=.2cm@C=.4cm{*={\bullet} \ar @{-} @/^{.7cm}/ [rrrrr]^{a}   & *={\bullet}  & *={\bullet} \ar @{=} @/^{.2cm}/ [r]^{e}  & *={\bullet} & *={\bullet} & *={\bullet}}}+\chi^{\ \xymatrix@R=.2cm@C=.4cm{*={\bullet} \ar @{-} @/^{.7cm}/ [rrrrr]^{a}   & *={\bullet}  \ar @{=} @/^{.3cm}/ [rr]^{e} & *={\bullet} & *={\bullet}  & *={\bullet} & *={\bullet}}}\bigg)\\
&\hspace*{6cm}\otimes
\bigg((2q-1)\chi^{\ \xymatrix@R=.2cm@C=.4cm{*={\bullet} & *={\bullet}  \ar @{-} @/^{.4cm}/ [rrr]^{b+d} & *={\bullet} & *={\bullet}  & *={\bullet} & *={\bullet}}}+(q-1)\sum_{f\in \FF_q^\times} \chi^{\ \xymatrix@R=.2cm@C=.4cm{*={\bullet} & *={\bullet}  \ar @{-} @/^{.6cm}/ [rrr]^{b}  & *={\bullet} \ar @{=} @/^{.2cm}/ [r]^{f} & *={\bullet}   & *={\bullet} & *={\bullet}}}\bigg)\\
&=
(2q-1)\chi^{\ \xymatrix@R=.2cm@C=.4cm{*={\bullet}  \ar @{-} @/^{.7cm}/ [rrrrr]^{a} & *={\bullet}  \ar @{-} @/^{.3cm}/ [rrr]^{b+d} & *={\bullet} & *={\bullet}  & *={\bullet} & *={\bullet}}}+(q-1)\sum_{f\in \FF_q^\times} \chi^{\ \xymatrix@R=.2cm@C=.4cm{*={\bullet}  \ar @{-} @/^{.9cm}/ [rrrrr]^{a} & *={\bullet}  \ar @{-} @/^{.6cm}/ [rrr]^{b}  & *={\bullet} \ar @{=} @/^{.2cm}/ [r]^{f} & *={\bullet}   & *={\bullet} & *={\bullet}}}+(2q-1)\sum_{e\in \FF_q^\times} \chi^{\ \xymatrix@R=.2cm@C=.4cm{*={\bullet} \ar @{-} @/^{.9cm}/ [rrrrr]^{a}   & *={\bullet}  \ar @{-} @/^{.5cm}/ [rrr]^{b+d} & *={\bullet} \ar @{=} @/^{.2cm}/ [r]^{e}  & *={\bullet} & *={\bullet} & *={\bullet}}}\\
&\hspace*{.5cm}+
(q-1)^2 \chi^{\ \xymatrix@R=.2cm@C=.4cm{*={\bullet} \ar @{-} @/^{.7cm}/ [rrrrr]^{a}   & *={\bullet}  \ar @{-} @/^{.4cm}/ [rrr]^{b} & *={\bullet}  & *={\bullet} & *={\bullet} & *={\bullet}}}+(q-1)(q-2)\sum_{g\in \FF_q^\times}\chi^{\ \xymatrix@R=.2cm@C=.4cm{*={\bullet} \ar @{-} @/^{.9cm}/ [rrrrr]^{a}   & *={\bullet}  \ar @{-} @/^{.6cm}/ [rrr]^{b} & *={\bullet} \ar @{-} @/^{.2cm}/ [r]^{g}  & *={\bullet} & *={\bullet} & *={\bullet}}}+(2q-1) \sum_{e\in \FF_q^\times}  \chi^{\ \xymatrix@R=.2cm@C=.4cm{*={\bullet} \ar @{-} @/^{.9cm}/ [rrrrr]^{a}   & *={\bullet}  \ar @{-} @/^{.5cm}/ [rrr]^{b+d} \ar @{=} @/^{.2cm}/ [rr]^(.7){e}  & *={\bullet}  & *={\bullet} & *={\bullet} & *={\bullet}}}\\
&\hspace*{.5cm}+
(q-1)\sum_{e,f\in \FF_q^\times}\chi^{\ \xymatrix@R=.2cm@C=.4cm{*={\bullet} \ar @{-} @/^{.9cm}/ [rrrrr]^{a}   & *={\bullet}  \ar @{-} @/^{.5cm}/ [rrr]^{b} \ar @{=} @/^{.2cm}/ [rr]^(.7){e}  & *={\bullet}  & *={\bullet} & *={\bullet} & *={\bullet}}}\\
&=
q(2q-1)\bigg(\chi^{\ \xymatrix@R=.2cm@C=.4cm{*={\bullet}  \ar @{-} @/^{.7cm}/ [rrrrr]^{a} & *={\bullet}  \ar @{-} @/^{.3cm}/ [rrr]^{b+d} & *={\bullet} & *={\bullet}  & *={\bullet} & *={\bullet}}}+\sum_{e\in \FF_q^\times} \chi^{\ \xymatrix@R=.2cm@C=.4cm{*={\bullet} \ar @{-} @/^{.9cm}/ [rrrrr]^{a}   & *={\bullet}  \ar @{-} @/^{.5cm}/ [rrr]^{b+d} & *={\bullet} \ar @{=} @/^{.2cm}/ [r]^{e}  & *={\bullet} & *={\bullet} & *={\bullet}}}\bigg)\\
&\hspace*{6cm}+
q(q-1)^2 \bigg(\chi^{\ \xymatrix@R=.2cm@C=.4cm{*={\bullet} \ar @{-} @/^{.7cm}/ [rrrrr]^{a}   & *={\bullet}  \ar @{-} @/^{.4cm}/ [rrr]^{b} & *={\bullet}  & *={\bullet} & *={\bullet} & *={\bullet}}}+\sum_{f\in \FF_q^\times} \chi^{\ \xymatrix@R=.2cm@C=.4cm{*={\bullet}  \ar @{-} @/^{.9cm}/ [rrrrr]^{a} & *={\bullet}  \ar @{-} @/^{.6cm}/ [rrr]^{b}  & *={\bullet} \ar @{=} @/^{.2cm}/ [r]^{f} & *={\bullet}   & *={\bullet} & *={\bullet}}}\bigg).
\end{align*} 

\noindent\textbf{Remark.}  The coefficients of the tensor products are not understood in general, although it is clear from Corollary \ref{TensorResult} that they are polynomial in $q$.

%
%

\subsection{Superinduction}

Let $S\subseteq \{1,2,\ldots,n\}$.  If $\mu\in \cS_S(q)$ and $\lambda\in \cS_n(q)$, then by Frobenius reciprocity,
$$\langle \chi^\lambda, \SInd_{U_S}^{U_n}(\chi^\mu)\rangle_{U_n}=\langle \Res_{U_S}^{U_n}(\chi^\lambda),\chi^\mu\rangle_{U_S}.$$
Thus, if 
\begin{equation*}
\SInd_{U_S}^{U_n}(\chi^\mu)  = \sum_{\nu} a_{\mu}^\nu \chi^\nu\quad\text{and}\qquad
 \Res_{U_S}^{U_n}(\chi^\lambda)  = \sum_{\gamma} b^\lambda_\gamma \chi^\gamma,
 \end{equation*}
 then by (\ref{OrthogonalityRelation})
 $$q^{|C(\lambda)|}a_{\mu}^\lambda=q^{|C(\mu)|}b^\lambda_\mu,$$
where $C(\nu)$ is the set of crossings of $\nu$.  Therefore,
$$\SInd_{U_S}^{U_n}(\chi^\mu)= \sum_{\nu} a_{\mu}^\nu \chi^\nu=\sum_{\nu} q^{|C(\mu)|-|C(\nu)|}b^\nu_\mu \chi^\nu.$$

In general, if $U_K\subseteq U_n$ with $K\in \cS_n$, then
\begin{equation}\label{SuperInductionByFrobenius}
\SInd_{U_K}^{U_n}(\chi^\mu)=\sum_{\nu} q^{|C_K(\mu)|-|C(\nu)|}b^\nu_\mu \chi^\nu,
\end{equation}
where $C_K(\nu)$ is the set of crossings that occur within the same parts of $K$.  

With this discussion, we obtain the following corollary of Sections \ref{SectionRestrictions} and \ref{SectionTensorProduct}.  When combined with Corollary \ref{SuperInductionTwoParts}, below, these results give a reasonably direct way to compute superinduction for some cases.

\begin{cor} \label{InducedFromTrivialTrivial}
Let $K=\{1,2,\ldots, k\}\cup\{k+1,k+2,\ldots,n\}\in \cS_n$ be a set-partition with two parts.  Then
$$\SInd_{U_K}^{U_n}(\One)=\sum_{{\lambda\in \cS_n(q)\atop \text{if $i\larc{}j\in \lambda$,}}\atop \text{then $i\nsim j\in K$}} q^{-|C(\lambda)|} \chi^\lambda,$$
where $i\sim j$ if and only if $i$ and $j$ are in the same part in $K$. 
\end{cor}

\begin{proof}
Note that by Theorem \ref{RestrictionRule} the coefficient of $\One=\chi^{\bullet\ \bullet\ \cdots\ \bullet}$ in $\Res_{U_K}^{U_n}(\chi^\lambda)$ is zero unless for every arc $i\larc{}j\in A(\lambda)$, the endpoints $i$ and $j$ are in different parts of $K$.   
Assume that every arc in $A(\lambda)$ satisfies this condition.  Then by Proposition \ref{LeviRestrictionRule}, 
\begin{align*}
\Res_{U_K}^{U_n}(\chi^\lambda)&=\bigg(\bigotimes_{i\larc{} k\in \lambda}\frac{\Res_{U_{K_1}}^{U_n}(\chi^\lambda)}{q^{|\{i<j<k\mid j\notin K_1\}|}} \bigg) \times \bigg(\bigotimes_{i\larc{} k\in \lambda}\frac{\Res_{U_{K_2}}^{U_n}(\chi^\lambda)}{q^{|\{i<j<k\mid j\notin K_2\}|}}\bigg) \\
&=\bigg(\bigotimes_{i\larc{} k\in \lambda} \One+\chi_{i,k}^{(1)}\bigg)\times\bigg(\bigotimes_{i\larc{} k\in \lambda} \One+\chi_{i,k}^{(2)}\bigg),
\end{align*}
where for $r\in \{1,2\}$,
$$\chi_{i,k}^{(r)}=\left\{\begin{array}{ll} \dd \sum_{i<j<k, j\in K_r\atop a\in \FF_q^\times} \chi^{j\larc{a}k}, & \text{if $k\in K_r$,}\\ 
 \dd \sum_{i<j<k, j\in K_r\atop a\in \FF_q^\times} \chi^{i\larc{a}j}, & \text{if $i\in K_r$,}\end{array}\right.$$
 
Thus, the coefficient of $\One$ in $\Res_{U_K}^{U_n}(\chi^\lambda)$ is $1$.  By (\ref{SuperInductionByFrobenius}),  the coefficient of $\chi^\lambda$ in $\SInd_{U_K}^{U_n}(\One)$ is $q^{|C_K(\bullet\ \bullet\ \cdots\ \bullet)|-|C(\lambda)|}=q^{-|C(\lambda)|}$.
\end{proof}

\noindent\textbf{Remarks.}   
\begin{enumerate}
\item[(a)] This corollary is a specific case.  The situation for even general two-part set partitions is more complicated (though seemingly tractable).  
\item[(b)]  Note that the superclass function $\SInd_{U_K}^{U_n}(\One)$ is a linear combination of supercharacters with fractional coefficients, so is a priori not generally a character. 
\end{enumerate}

Corollary \ref{InducedFromTrivialTrivial} has some immediate combinatorial consequences.  Let
$$SG_{n\times m}=\{a\in M_{n\times n}(\{0,1\})\mid \text{$a$ has at most one 1 and every row and column}\}$$
be the set of $m\times n$ $0$-$1$ matrices with at most one 1 in every row an column.  Define statistics for $w\in SG_{m\times n}$
\begin{align*}
\mathrm{ones}(w) &= |\{(i,j)\in [1,n]\times [1,m]\mid w_{ij}=1\}|\\
\mathrm{sow}(w) &= |\{(j,k)\in [1,n]\times [1,m]\mid w_{jk}=0, \text{$w_{ik}=1$ for some $i<j$ or $w_{jl}=1$ for some $k<l$}\}|.
\end{align*}
For example, if 
$$w=\left(\begin{array}{cccc} \underline{0} & 1 & 0 & 0 \\ \underline{0} & \underline{0} & \underline{0} & 1\\ 0 & \underline{0} & 0 & \underline{0} \end{array}\right), \qquad \text{then} \qquad \begin{array}{l}\mathrm{ones}(w)=2\\
\mathrm{sow}(w)=6. \end{array}$$

\begin{cor}  Let $m$ and $n$ be positive integers.  Then
\begin{enumerate}
\item[(a)] $\dd q^{mn}=\sum_{w\in SG_{m\times n}} (q-1)^{\mathrm{ones}(w)}q^{\mathrm{sow}(w)}$
\item[(b)] $\dd 0 = \sum_{w\in SG_{m\times n}} (-1)^{w_{1n}} (q-1)^{\mathrm{ones}(w)}q^{\mathrm{sow}(w)}.$
\end{enumerate}
\end{cor}
\begin{proof}
(a) We evaluate the equation in Corollary \ref{InducedFromTrivialTrivial} at the identity of $U_{m+n}$, when $K=\{1,2,\ldots, m\} \cup\{m+1,\ldots, m+n\}$.  The degree of $\SInd_{U_K}^{U_n}(\One)$ is the size of the factor group $|U_{m+n}/U_K|=q^{mn}$.  

On the right-hand side of the equation we sum over the set $X$ of all possible set-partitions whose arcs $i\larc{}l$ satisfy $i\in \{1,\ldots, m\}$ and $l\in \{m+1,\ldots, m+n\}$.  Consider the surjective map $\varphi:X \rightarrow SG_{m\times n}$ given by
$$\varphi(\lambda)_{i,l-m} =1 \qquad \text{if and only if}\qquad i\larc{}l\in A(\lambda).$$
Note that the size of the preimage of $w\in SG_{m\times n}$ is $|\varphi^{-1}(w)|=(q-1)^{\mathrm{ones}(w)}$.  For $\lambda\in X$,
$$\chi^\lambda(1)=\prod_{i\larc{}l\in A(\lambda)} q^{l-i-1}=\prod_{i\larc{}l\in A(\lambda)} q^{l-m-1+m-i}=q^{\mathrm{sow}(\varphi(\lambda))} q^{|C(\lambda)|},$$
where each crossing in $C(\lambda)$ accounts for an over-counting of an entry in the computation of $\mathrm{sow}(\varphi(\lambda))$.  Thus,
$$q^{mn}=\SInd_{U_K}^{U_n}(\One)(1)=\sum_{{\lambda\in \cS_n(q)\atop \text{if $i\larc{}j\in \lambda$,}}\atop \text{then $i\nsim j\in K$}} q^{-|C(\lambda)|} \chi^\lambda(1)=\sum_{w\in SG_{m\times n}} (q-1)^{\mathrm{ones}(w)}q^{\mathrm{sow}(w)}.$$

(b) Instead of evaluating at the identity, as in (a), evaluate both sides at one of the central non-identity superclasses, each of whose set-partition $\mu$ is of the form $\mu=1\larc{a}(m+n)$ for some $a\in \FF_q^\times$.   
\end{proof}

We conclude with some observations relating superinduction to these superpermutation ``characters."  

\begin{proposition}
Let $H\subseteq G$ be pattern groups.  For a superclass functions $\gamma$ and $\eta$ of $G$,
$$\langle \gamma \otimes \SInd_H^G(\One),\eta\rangle=\langle \Res_H^G(\gamma),\Res_H^G(\eta)\rangle.$$
\end{proposition}
\begin{proof}
By a typical Frobenius reciprocity argument,
\begin{align*}
\langle \gamma\otimes \SInd_H^G(\One),\eta\rangle_G & = \frac{1}{|G|}\sum_{g\in G}\gamma(g)\SInd_H^G(\One)(g)\overline{\eta(g)}\\
&= \frac{1}{|G|^2|H|}\sum_{g,x,y\in G} \gamma(g)\overset{\circ}{\One}(x(g-1)y+1)\overline{\eta(g)}\\
\intertext{where $\overset{\circ}{\One}(g)$ is 1 if $g\in H$ and $0$ otherwise.  By reindexing and using the fact that $\gamma$ and $\eta$ are superclass functions of $G$,}
&= \frac{1}{|G|^2|H|}\sum_{k,x,y\in G} \gamma(x^{-1}(k-1)y^{-1}+1)\overset{\circ}{\One}(k)\overline{\eta(x^{-1}(k-1)y^{-1}+1)}\\
&=\frac{1}{|G|^2|H|}\sum_{k,x,y\in G} \gamma(k)\overset{\circ}{\One}(k)\overline{\eta(k)}\\
&=\frac{1}{|H|}\sum_{k\in H} \gamma(k)\overline{\eta(k)}\\
&=\langle \Res_H^G(\gamma),\Res_H^G(\eta)\rangle_H.\qedhere
\end{align*}
\end{proof}

By choosing $\gamma$ in the above proposition appropriately, we get the following useful corollary via Frobenius reciprocity.

\begin{cor}\label{SuperinductionPermutationCharacter}
Let $H\subseteq G$ be pattern groups, and let $\mu\in (H-1)$.  If $\chi^\mu(1)\Sinf_H^G(\chi^\mu)(h)=\Sinf_H^G(\chi^\mu)(1)\chi^\mu(h)$, for all $h\in H$, then
$$\SInd_H^G(\chi^\mu)=\frac{\chi^\mu(1)}{\Sinf_H^G(\chi^\mu)(1)}\Sinf_H^G(\chi^\mu)\otimes \SInd_H^G(\One).$$
\end{cor}

The assumption in Corollary \ref{SuperinductionPermutationCharacter} is not so unusual.  In fact,

\begin{cor} \label{SuperInductionTwoParts}
Let $U_K\subseteq U_L$ be parabolic subgroups of $U_n$, where $K=K_1\cup K_2\cup\cdots\cup K_\ell,L\in \cS_n$.  Then for $\mu\in \cS_{K_1}(q)\times \cS_{K_2}(q)\times\cdots\times \cS_{K_\ell}(q)$,
$$\SInd_{U_K}^{U_L}(\chi^\mu)=\frac{\chi^\mu(1)}{\Sinf_{U_K}^{U_L}\big(\chi^{\mu}\big)(1)}\Sinf_{U_K}^{U_L}(\chi^\mu)\otimes \SInd_{U_K}^{U_L}(\One).$$
\end{cor}
\begin{proof}
Note that for $i<j\in \cP_K$ and $k\in U_K\subseteq U_L$,
$$\chi^{\langle i\larc{a}j\rangle_L}(k)=\frac{\chi^{\langle i\larc{a}j\rangle_L}(1)}{\chi^{\langle i\larc{a}j\rangle_K}(1)}\chi^{\langle i\larc{a}j\rangle_K}(k).$$
Thus, by the decomposition of characters, for $\mu\in \cS_{K_1}(q)\times \cS_{K_2}(q)\times\cdots\times \cS_{K_\ell}(q)$  and $k\in U_K$,
$$\Sinf_{U_K}^{U_L}(\chi^\mu)(k)=\chi^{\langle\mu\rangle_L}(k)=\frac{\chi^{\langle\mu\rangle_L}(1)}{\chi^{\mu}(1)}\chi^{\mu}(k),$$
and the result now follows from Corollary \ref{SuperinductionPermutationCharacter}.
\end{proof}

\noindent\textbf{Remark.}  While the assumption in Corollary \ref{SuperinductionPermutationCharacter} is sufficient, it is not necessary.  For example, if
$$H=\left\{\left(\begin{array}{ccc} 1 & \ast & \ast\\ 0 & 1 & 0 \\ 0 & 0 & 1\end{array}\right)\right\}\subseteq U_3=\left\{\left(\begin{array}{ccc} 1 & \ast & \ast\\ 0 & 1 & \ast \\ 0 & 0 & 1\end{array}\right)\right\}.$$
then for these groups, 
$$\chi^{\left(\begin{smallmatrix} 0 & 0 & 1\\ 0 & 0 & 0 \\ 0 & 0 & 0\end{smallmatrix}\right)}\left(\begin{array}{ccc} 1 & 1 & 0\\ 0 & 1 & 0 \\ 0 & 0 & 1\end{array}\right)=1
\qquad \text{and}\qquad 
\Sinf_H^{U_3}\bigg(\chi^{\left(\begin{smallmatrix} 0 & 0 & 1\\ 0 & 0 & 0 \\ 0 & 0 & 0\end{smallmatrix}\right)}\bigg)\left(\begin{array}{ccc} 1 & 1 & 0\\ 0 & 1 & 0 \\ 0 & 0 & 1\end{array}\right)=0.$$
However, it remains true that
$$\SInd_H^G\bigg(\chi^{\left(\begin{smallmatrix} 0 & 0 & 1\\ 0 & 0 & 0 \\ 0 & 0 & 0\end{smallmatrix}\right)}\bigg)=q^{-1}\Sinf_H^G\bigg(\chi^{\left(\begin{smallmatrix} 0 & 0 & 1\\ 0 & 0 & 0 \\ 0 & 0 & 0\end{smallmatrix}\right)}\bigg)\otimes \SInd_H^G(\One).$$
In fact,  the conclusion of Corollary \ref{SuperinductionPermutationCharacter} may be true for all pattern groups; I know of no counter-example.

%
%

\end{document}